\documentclass[11pt]{iopart}
\usepackage{iopams}
\usepackage{amsfonts,amsbsy,amssymb}
\usepackage{amssymb}
\usepackage{amsthm}
\usepackage{graphicx}
\usepackage[format=plain, textfont=it]{caption}
\usepackage{subfig}
\usepackage{pgfplotstable}
\usepackage{tikz}

\input{tcilatex}
\begin{document}

\title{Single measurement experimental data for an inverse medium problem
inverted by a multi-frequency globally convergent numerical method}
\author{Aleksandr E. Kolesov$^{1,2}$, Michael V. Klibanov$^1$, Loc H. Nguyen$%
^1$, Dinh-Liem Nguyen$^1$ and Nguyen T. Th\`{a}nh$^3$}

\begin{abstract}
The recently developed globally convergent numerical method for an inverse
medium problem for the Helmholtz equation proposed in \cite{KLN} is tested
on experimental data. The data were originally collected in the time domain,
whereas the method works in the frequency domain with the multi-frequency
data. Due to a huge discrepancy between the collected and computationally
simulated data, the straightforward Fourier transform of the experimental
data does not work. Hence, it is necessary to develop a heuristic data
preprocessing procedure. This procedure is described. The preprocessed data
are used as the input for the inversion algorithm. Numerical results
demonstrate good accuracy in the reconstruction of both refracive indices
and locations of targets. Furthermore, the reconstruction errors for
refractive indices of dielectric targets are significantly less than errors
of a posteriori direct measurements.
\end{abstract}

\address{$^1$ Department of Mathematics and Statistics, University of North
Carolina at Charlotte, Charlotte, NC 28223, USA} 
\address{$^2$ Institute of
Mathematics and Information Science, North-Eastern Federal University,
Yakutsk, Russia} 
\address{$^3$ Department of Mathematics, Iowa State
University, Ames, IA 50011, USA} 
\ead{akolesov@uncc.edu, mklibanv@uncc.edu,
lnguyen50@uncc.edu, dnguyen70@uncc.edu, thanh@iastate.edu}

\pagestyle{plain}

\vspace{10pt} 

%
%
%
%
%

\textbf{Key Words:} experimental time dependent data, multi-frequency data,
global convergence, coefficient inverse problems, inverse scattering problems

\textbf{2010 Mathematics Subject Classification:} 35R30, 78A46

\section{Introduction}

\label{sec:1}

In this work we consider an inverse medium scattering problem for the
Helmholtz equation in the three dimensional space $\mathbb{R}^{3}$. The
objective is to reconstruct the coefficient of the Helmholtz equation in a
bounded domain. The coefficient represents the spatially distributed
dielectric constant of the medium. Our target application is in the
detection and identification of explosives, such as antipersonnel mines and
improvised explosive devices (IEDs). We calculate dielectric constants and
estimate locations of objects which mimic explosives. 
Guided by our target application, we use only a single boundary measurement
of the backscatter wave.

Currently, the radar community relies mainly on the intensity of the radar
images, which are obtained by migration-type imaging methods, to obtain
geometrical information such as the shapes, the sizes, and the locations of
the targets, see, e.g., \cite{KSNF1, Soumekh:1999,Yilmaz:1987}. Hence, the
additional information about values of dielectric constants of targets of
interest might help in the future to develop classification algorithms,
which would better differentiate between explosives and clutter. The targets
in our experiments are located in air. It is known that, for example IEDs
can be located in air. On the other hand, the case of targets buried in the
ground is one of goals of our future research. We also note that this case
was studied in \cite{TBKF2} for experimental time dependent data using the
previously developed globally convergent inverse algorithm of \cite{BK,BK1}.

Another term for the inverse medium problem under consideration is
Coefficient Inverse Problem (CIP). Recently a globally convergent numerical
method for this CIP with multi-frequency data resulting from a single
measurement event has been developed by this group in \cite{KLN}. This
method was tested on computationally simulated data in \cite{KLN} and on
experimental multi-frequency data in \cite{Liem}. Obviously, testing of an
inversion algorithm on several types of experimental data is a good idea
since it provides better assurances of the performance of this algorithm.
Thus, the goal of the current paper is to test the technique of \cite{KLN}
on time dependent experimental data. These data were collected on a
microwave scattering facility in the University of North Carolina at
Charlotte and were used in \cite{IPexp1,TBKF1} to test a different globally
convergent inverse algorithm of \cite{BK,BK1}. 
Note that, motivated by our target application, we measured only the
backscatter data generated by a source located at a fixed location. Thus,
this is a single measurement data, which is one of the most challenging
cases for any inversion algorithm.

It seems to be at the first glance that the easiest way to apply the
frequency domain globally convergent algorithm of \cite{KLN} to time domain
data is to apply the Fourier transform to the data and then to use the
resulting data as the input for the algorithm. However, the straight forward
application of this idea does not work here. The latter is due to a 
\underline{huge discrepancy} between experimental and computationally
simulated data, see section 4.2.1. Such a discrepancy was observed earlier
for both the time dependent \cite{IPexp1,TBKF1,TBKF2} and multi-frequency 
\cite{Liem} experimental data. We remark here that conventional data
denoising techniques do not work for our data because of its complicated
structure, see section 4.2.1. Therefore, it was concluded in \cite{IPexp1,
Liem, TBKF1, TBKF2} that a heuristic data preprocessing procedure is
necessary to make the preprocessed data look at least somehow similar to the
computationally simulated ones. The data preprocessing procedure of the
current paper is described in section 4. The result of this procedure is
used as the input for the algorithm of \cite{KLN}.

The first step of the globally convergent numerical method of \cite%
{BK,BK1,IPexp1,TBKF1,TBKF2} is the application of the Laplace transform with
respect to time to the solution of a hyperbolic wave-like PDE. However, it
was observed in \cite{Liem} that this technique does not work for the
multi-frequency experimental data of \cite{Liem}. The reason of this is that
these data are stable only on a small interval of frequencies concentrated
around a certain \textquotedblleft optimal" frequency. A similar observation
is made in section 4.2.7 for the Fourier transform of our preprocessed time
dependent data. We point out that this phenomenon is not observed in
computationally simulated data. Therefore, this is one of significant
discrepancies between real and simulated data. In other words, when working
with our data, one can rely only on a small interval of frequencies. We note
that the latter is one of conditions of the global convergence theorem of 
\cite{KLN}. The integral of the inverse Fourier transform, which is carried
out over only that small interval of frequencies, cannot provide a
reasonable accuracy of resulting time dependent data.

A different type of experimental data for multi-frequencies and multiple
locations of sources was collected in Fresnel Institute (Marseille, France).
These data were used then by several teams to solve CIPs, see a summary of
results in \cite{Litman}. In this case an anechoic chamber was used for the
data collection. The latter removed many parasitic signals from the data. As
a result, these data matched well with simulations \cite{Geff}. However,
keeping in mind our target application, we have collected the data in a more
realistic environment of a regular office room filled with the office
furniture, computers, etc. Naturally, this resulted in reflections of our
signals from various items of this room. The latter, in turn led to the
above mentioned huge discrepancy.

The numerical methods of \cite{BK,BK1} and \cite{KLN} are the so-called
\textquotedblleft approximately globally convergent methods". A detailed
discussion of the notion of the approximate global convergence can be found
in \cite{BK,BK1}. We now explain this notion briefly. Any CIP is a highly
nonlinear problem. Therefore, an important question in its numerical
treatment is: \emph{How to obtain at least one point in a sufficiently small
neighborhood of the exact solution without any advanced knowledge of this
neighborhood}? A numerical method for a CIP is called \textquotedblleft
approximately globally convergent" (globally convergent in short, or GCM) if
a theorem is proved, which claims that, under a certain reasonable
mathematical assumption, this method addresses the above question
positively, i.e. it delivers that point. We call this theorem the
\textquotedblleft global convergence theorem". The estimate of the distance
between that point and the true coefficient should depend on the error in
the data and some parameters of the discretization. We point out that the
proximity of that point to the true coefficient is the main advantage of the
GCM. Indeed, as soon as such a point is found, the solution can be refined
via a small perturbation approach, see, e.g. Chapter 4 of \cite{BK}. That
notion of a reasonable mathematical assumption is well justified by the well
known fact that the goal of the development of such numerical methods for
CIPs, which would positively address the above question, is a \underline{%
tremendously challenging} one, especially for the case of a single
measurement data. We refer to \cite{BK,BK1} for detailed discussions of the
notion of the approximate global convergence.

We note that there is a vast literature on reconstruction methods developed
for solving CIPs. To study the CIP under consideration, in which weak
scattering assumptions are not applicable, the probably best known approach
is nonlinear optimization schemes, see, e.g.,~\cite{Chave2010} and
references therein. However, it is well-known that the methods based on
nonlinear optimization schemes heavily rely on a strong \emph{a priori}
knowledge about the target. In particular, the convergence of those methods
requires a good \emph{a priori} initial approximation of the exact solution,
that is, the starting point of iterations should be chosen to be
sufficiently close to the solution. Hence, we call such methods
\textquotedblleft locally convergent". Note that in our desired applications
such \emph{a priori} knowledge is not always available. Concerning
qualitative reconstruction methods for inverse scattering problems, we refer
to~\cite{Amma2016, Cakon2006, Dorn2015, Kirsc2008, Li2014, Potth2006,
Griesmaier:IP2011, Potthast:IP2010} and the references therein. These
methods do not require good first guesses. However, they reconstruct only
the shapes of scattering objects instead of their material properties.
Another approach to inverse scattering problems using multi-frequency data
is so called frequency-hopping algorithms. These methods start from a low
frequency to get a rough reconstruction, then improve it recursively using
higher frequencies. By using a low frequency at the first iteration, it is
not necessary to start from a good first guess. However, their convergence
is known only in some limited cases. We refer the reader to \cite{Bao2015,
Chen:IP1997, C-L:IEEE1995, Sini2012, ST:ESAIMM2NA2014} for these methods.

%
%




Finally, we refer to different globally convergent numerical methods for
solving CIPs with multiple measurements of the Dirichlet-to-Neumann map \cite%
{Agalt2015, Agalt2014, Belis2016, Ivano2016, Kaban2015, Kaban2004, Novik2015}%
. These techniques were tested on computationally simulated data in~\cite%
{Belis2016, Kaban2015,Kaban2004}. We recall that our GCM in this paper deals
with only a single backscatter measurement.

The paper is organized as follows. We present in section 2 the formulation
of the direct and inverse problems considered in this paper. Section 3 is a
brief summary of our GCM. The data collection and preprocessing are
described in section 4. We present in section 5 the numerical implementation
of the method for the preprocessed data. Section 6 is devoted to numerical
results obtained from the numerical implementation in section 5. Discussion
is present in section 7.

\section{Problem statement}

\label{sec:2}

\subsection{Forward and inverse problems}

\label{sec:2.1} We state in this section the forward and inverse problems
under the consideration. We consider dimensionless variables in the sections
involving the theory of the globally convergent method (sections 2 and 3).
However, since the aim of the paper is to invert the measured data, which
have dimensions, we will explain in section 4.1 how do we make variables
dimensionless.

Let $x=\left( x_{1},x_{2},x_{3}\right) \in \mathbb{R}^{3}.$ Let $%
B(R)=\{|x|<R\}\subset \mathbb{R}^{3}$ be the ball of the radius $R$ with the
center at $\left\{ x=0\right\} $ and $\Omega _{1}\Subset \Omega \Subset B(R)$
be bounded domains. We assume that the medium is isotropic and non magnetic.
Let the function $c(x)$ be the smooth spatially distributed dielectric
constant satisfying the following conditions: 
\begin{equation}
c(x)\geq 1,x\in \mathbb{R}^{3},\quad \mbox{and }c(x)=1,x\in \mathbb{R}%
^{3}\setminus \Omega _{1}.  \label{eq:beta}
\end{equation}%
%
%
Condition (\ref{eq:beta}) means that the dielectric constant of the
scattering medium $\Omega _{1}$ exceeds the one of air (where $c(x)\approx
1) $ and that outside of the domain $\Omega $ $\ $we have air. The
smoothness condition of $c(x)$ is a technical one.\ It was imposed in \cite%
{KR,KLN} to justify the asymptotic behavior of the solution of the Helmholtz
equation~(\ref{eq:helmholtz}).

It was pointed out in Chapter 13 of the classical book \cite{Born} that if
the dielectric constant varies sufficiently slowly on the scale of the
wavelength, then the solution of the Maxwell's equations can be well
approximated by the solution of the scalar Helmholtz equation for a certain
component of the electric field $E=(E_{x_{1}},E_{x_{2}},E_{x_{3}})$. Hence,
we work below only with the single Helmholtz equation. In our experiments
the incident electric wave field was $E^{inc}=(0,E_{x_{2}}^{inc},0)$ and
only the backscatter component of $E_{x_{2}}\left( x,t\right) $ was
measured. Hence, in the frequency domain we set $u\left( x,k\right) =%
\widetilde{E}_{x_{2}}\left( x,k\right) $, where $\widetilde{E}_{x_{2}}\left(
x,k\right) $ is the frequency domain analog of $E_{x_{2}}\left( x,t\right) .$
The forward problem for the Helmholtz equation is 
\begin{equation}
\Delta u+k^{2}c(x)u=0,\quad x\in \mathbb{R}^{3},  \label{eq:helmholtz}
\end{equation}%
where $k=2\pi /\lambda $ is the wave number and $\lambda $ is the
dimensionless wavelength. Let 
\begin{equation}
u_{0}\left( x,k\right) =e^{ikx_{3}}  \label{3}
\end{equation}%
be the incident plane wave propagating along the $x_{3}-$axis. Let $%
u_{sc}(x,k)$ be the scattered wave due to the heterogeneous scattering
medium $\Omega _{1}.$ We seek the solution of equation (\ref{eq:helmholtz}),
i.e. the full wave field, in the following form: 
\begin{equation}
u(x,k)=u_{0}(x,k)+u_{sc}(x,k),  \label{eq:total}
\end{equation}%
\begin{equation}
\frac{\partial u_{sc}}{\partial r}-iku_{sc}=o(r^{-1}),\quad r=|x|\rightarrow
\infty ,  \label{eq:radiation}
\end{equation}%
where (\ref{eq:radiation}) is the Sommerfeld radiation condition. It is well
known that for each $k>0$ there exists unique solution $u\in C^{2}\left( 
\mathbb{R}^{3}\right) $ of the problem (\ref{eq:helmholtz})-(\ref%
{eq:radiation}).

The problem (\ref{eq:helmholtz})-(\ref{eq:radiation}) is our forward
problem. We now pose the inverse problem. Let $\underline{k}$ and $\overline{%
k}$ be two positive constants and $\underline{k}<\overline{k}$. Let $\Gamma
\subset \partial \Omega $ be the backscatter part of the boundary of the
domain $\Omega .$ Our CIP is stated as:

\textbf{CIP for multi-frequency data}. \textit{Determine the coefficient }$%
c(x)$\textit{\ for }$x\in \Omega _{1}$\textit{, given the backscatter data }$%
g(x,k)$ \emph{on} $\Gamma ,$\textit{\ } 
\begin{equation}
g(x,k)=u(x,k),\quad x\in \Gamma ,\,k\in \lbrack \underline{k},\overline{k}].
\label{eq:cip_fdd}
\end{equation}

Along with this CIP, it is convenient to consider its time domain analog. In
this case the forward problem is 
\begin{equation}
c(x)w_{tt}=\Delta w,\quad (x,t)\in \mathbb{R}^{3}\times (-\infty ,\infty ),
\label{4}
\end{equation}%
\begin{equation}
w\left( x,t\right) =\delta \left( t-x_{3}\right) +\widetilde{w}\left(
x,t\right) ,  \label{5}
\end{equation}%
\begin{equation}
\widetilde{w}\mid _{t<-R}=0.  \label{6}
\end{equation}%
In (\ref{4})-(\ref{6}) $w\left( x,t\right) $ is the total wave field, $%
\delta \left( t-x_{3}\right) $ is the plane wave propagating along the $%
x_{3}-$axis and $\widetilde{w}$ is the scattered wave field. Condition (\ref%
{6}) means that the total wave field $w\left( x,t\right) =\delta \left(
t-x_{3}\right) $ for those times $t$ which are prior the moment of time when
this plane wave reaches the scattering medium $\Omega _{1}$. Results of \cite%
{V} ensure that, under certain conditions, the function $w\left( x,t\right) $
decays exponentially as $t\rightarrow \infty ,$ together with its
appropriate derivatives, uniformly for all points $x$ belonging to an
arbitrary selected bounded domain $G\subset \mathbb{R}^{3}.$ Using this fact
as well as the smoothness of $c(x)$, it was proven in \cite{KR} that\textbf{%
\ }the functions $u$ and $w$ are connected via the Fourier transform,

\begin{equation}
u(x,k)=\displaystyle\int\limits_{-\infty }^{\infty }w(x,t)\,e^{ikt}\,dt.
\label{7}
\end{equation}%
The latter likely works well to transform the computationally simulated data
for the problem (\ref{4})-(\ref{6}) in the frequency domain and to apply the
method of \cite{KLN} then. However, before applying (\ref{7}) to our
experimental data, we have to preprocess them, see sections 4.2.1-4.2.5.

As to the uniqueness of the above CIP, this is a CIP with a single
measurement data. Uniqueness theorems for multidimensional CIPs with single
measurement data are currently proven only via the method, which was
originated in \cite{BukhKlib}. This method is based on Carleman estimates.
There are many publications on this technique of a number of authors about
uniqueness theorems for CIPs with a finite number of measurements, see, e.g. 
\cite{Ksurvey,Yam} for surveys. However, if applying this method to our CIP,
then one needs to assume that zero in the right hand side of equation (\ref%
{eq:helmholtz}) is replaced with a function $f\left( x\right) $ which is non
zero everywhere in $\overline{\Omega }_{1}.$ In fact, it is well known that
it is a very challenging open problem: to prove uniqueness of our CIP when
zero is in the right hand side of (\ref{eq:helmholtz}). Hence, we assume
below that the uniqueness holds for our CIP.

\subsection{The Lippmann-Schwinger equation}

\label{sec:2.2}

In our inverse algorithm we need to solve the forward problem (\ref%
{eq:helmholtz})-(\ref{eq:radiation}) for $k=\overline{k},$ with different
coefficients $c\left( x\right)$ obtained in the iterative process. We are
doing this via solving the Lippmann-Schwinger equation 
\begin{equation}
u(x,\overline{k})=\exp (i\overline{k}x_{3})+\overline{k}^{2}\int_{\Omega
}\Phi (x,y,\overline{k})\left( c\left( y\right) -1\right) u(y,\overline{k}%
)dy,\quad x\in \Omega .  \label{eq:lippmann-schwinger}
\end{equation}%
Here $\Phi (x,y)$ is the fundamental solution of the Helmholtz equation with 
$c\left( x\right) \equiv 1$ for $k=\overline{k},$ 
\[
\Phi (x,y,\overline{k})=\frac{\exp (i\overline{k}|x-y|)}{4\pi |x-y|},\quad
x\neq y. 
\]
It is well known that if the function $c$ is smooth and satisfies condition (%
\ref{eq:beta}), then equation (\ref{eq:lippmann-schwinger}) is equivalent to
the problem (\ref{eq:helmholtz})-(\ref{eq:radiation}), see Chapter 8 in \cite%
{CK}. To solve equation (\ref{eq:lippmann-schwinger}), we use a fast
numerical method developed in \cite{Vainikko}.

\section{Globally convergent numerical method}

\label{sec:3}

In this section, we refer to \cite{KLN} for all theoretical details. Even
though the measured data were collected only on the backscatter part $\Gamma 
$ of the boundary $\partial \Omega ,$ the theory of \cite{KLN} works only
for the case when the data are given on the entire boundary $\partial \Omega
.$ Thus, we assume the latter in the current section. We show in section 5
how to make this method work for backscatter data.

\subsection{Nonlinear integro-differential equation}

\label{sec:3.1}

The first step of the GCM consists of obtaining a nonlinear
integro-differential equation. This step is not a part of locally convergent
algorithms. In fact, this idea is actually taken from the method of \cite%
{BukhKlib}. It was established in \cite{KLN} that for sufficiently large $k$
the function $u(x,k)\neq 0$ for all $x\in \overline{\Omega }.$ Hence, we
assume from now on that numbers $\underline{k},\overline{k}$ are
sufficiently large and that $k\in \left[ \underline{k},\overline{k}\right] .$
In \cite{KLN} the unique function $v\left( x,k\right) \in C^{2}\left( 
\overline{\Omega }\right) $ was constructed such that 
\begin{equation}
u(x,k)=e^{v\left( x,k\right) },x\in \overline{\Omega }.  \label{eq:v}
\end{equation}%
Substituting (\ref{eq:v}) in equation (\ref{eq:helmholtz}), we obtain 
\begin{equation}
\Delta v+(\nabla v)^{2}=-k^{2}c\left( x\right) ,x\in \overline{\Omega }.
\label{eq:c_from_v}
\end{equation}%
Introduce the function $q\left( x,k\right) $ as%
\begin{equation}
q(x,k)=\partial _{k}v(x,k)=\frac{\partial _{k}u(x,k)}{u(x,k)}.  \label{eq:q}
\end{equation}%
Then, differentiating equation (\ref{eq:c_from_v}) with respect to $k$, we
obtain 
\begin{equation}
\Delta q+2\nabla q\nabla v=\frac{2(\Delta v+(\nabla v)^{2})}{k}.
\label{eq:q_v_eq}
\end{equation}%
By (\ref{eq:q}) 
\begin{equation}
v(x,k)=-\int_{k}^{\overline{k}}q(x,\kappa )d\kappa +V(x),  \label{eq:v_q}
\end{equation}%
\begin{equation}
V(x)=v(x,\overline{k}).  \label{8}
\end{equation}%
The function $V(x)$ is called the \textquotedblleft tail function" and the
truncation frequency $\overline{k}$ plays the role of a regularization
parameter.

Substituting (\ref{eq:v_q}) in (\ref{eq:q_v_eq}), we obtain the following
nonlinear integro-differential equation for the function $q$ 
\[
k\Delta q+2\nabla q\left( -\int_{k}^{\overline{k}}\nabla q(x,\kappa )d\kappa
+\nabla V\right) = 
\]%
\begin{equation}
2\left( -\int_{k}^{\overline{k}}\Delta q(x,\kappa )d\kappa +\Delta V+\left(
-\int_{k}^{\overline{k}}\nabla q(x,\kappa )d\kappa +\nabla V\right)
^{2}\right).  \label{9}
\end{equation}%
Equation (\ref{9}) is complemented with the Dirichlet boundary condition 
\begin{equation}
q(x,k)=\frac{\partial _{k}g(x,k)}{g(x,k)}=:\psi (x,k),\quad x\in \partial
\Omega .  \label{10}
\end{equation}

In (\ref{9}), (\ref{10}) both functions $q$ and $V$ are unknown and need to
be approximated. We approximate them via a predictor-corrector procedure, in
which approximations for $V$ are predictors and approximations for $q$ are
correctors. As the zero step, we find a good approximation $V_{0}$ for the
function $V$. Then an iterative process is used to update $q,c,u$ and $V$.
More precisely, given an approximation of the tail function $V$, we solve
the Dirichlet boundary value problem for an elliptic equation in the domain $%
\Omega $ to find the next approximation for the function $q$. Next, we
update the approximation for the function $v$ via (\ref{eq:v_q}) and then
for the unknown coefficient $c\left( x\right) $ via (\ref{eq:q_v_eq}). Next,
given that update for $c$, we solve the Lippman-Schwinger equation (\ref%
{eq:lippmann-schwinger}) for $u(x,\overline{k})$. Observe that (\ref{9})
implies that we do not need to know the function $V$ itself. Rather, we need
to know its gradient and its Laplacian, which are given by 
\begin{equation}
\nabla V\left( x\right) =\frac{\nabla u\left( x,\overline{k}\right) }{%
u\left( x,\overline{k}\right) },\Delta V\left( x\right) =\func{div}\left( 
\frac{\nabla u\left( x,\overline{k}\right) }{u\left( x,\overline{k}\right) }%
\right) .  \label{200}
\end{equation}

\subsection{Discretization with respect to $k$}

\label{sec:3.2}

To approximate both functions \ $q$ and $V$ in the above outlined iterative
process, we use a discretization of equation (\ref{8}) with respect to $k\in
\lbrack \underline{k},\overline{k}].$ We divide the interval $[\underline{k},%
\overline{k}]$ into $N$ uniform subintervals with the discretization step
size $h$, 
\begin{equation}
\underline{k}=k_{N}<k_{N-1}<\dots <k_{1}<k_{0}=\overline{k},\quad
h=k_{n-1}-k_{n}.  \label{11}
\end{equation}%
Next, we approximate $q(x,k)$ by a piecewise constant function with respect
to $k\in \lbrack \underline{k},\overline{k}]$. This and (\ref{10}) imply
that the function $\psi (x,k)$ should also be approximated by a piecewise
constant function with respect to $k$. In other words, 
\begin{equation}
q(x,k) \simeq q_{n}(x),\quad \psi (x,k) \simeq \psi _{n}(x),\quad k\in
\lbrack k_{n},k_{n-1}),\quad n=1,\dots N.  \label{eq:q_psi_appr}
\end{equation}%
We set $q_{0}(x)\equiv 0$ and denote 
\[
\overline{q_{n-1}}=\sum_{j=0}^{n-1}q_{j}(x). 
\]%
Hence, 
\[
\int_{k}^{\overline{k}}q(x,\kappa )d\kappa \simeq (k_{n-1}-k)q_{n}+h%
\overline{q_{n-1}},\quad k\in \lbrack k_{n},k_{n-1}). 
\]%
Now we set in (\ref{9}) $k\in \left( k_{n},k_{n-1}\right] $ and assume that $%
h$ is sufficiently small. In particular, we assume that $h\overline{k}\ll 1$%
. Hence, we eliminate those terms in (\ref{9}), whose absolute values are $%
O(h)$ as $h\rightarrow 0$. Next, to eliminate the dependence on the varying
parameter $k$, we integrate the resulting equation with respect to $k\in
\left( k_{n},k_{n-1}\right] $ and then divide by $h$. We obtain after some
manipulations 
\begin{eqnarray}  \label{12}
&\Delta q_{n}-A_{n}h\nabla \overline{q_{n-1}}\nabla q_{n}=Q_{n},\quad x\in
\Omega,  \nonumber \\
&Q_{n}=-A_{n}\nabla q_{n-1}\nabla V_{n-1}+\frac{2(\Delta V_{n-1}+(\nabla
V_{n-1})^{2})}{k_{n-1}}  \nonumber \\
& \qquad \, -\frac{4\nabla V_{n-1}h\nabla \overline{q_{n-1}}}{k_{n-1}}-\frac{%
2h\Delta \overline{q_{n-1}}}{k_{n-1}}, \\
&q_{n}=\psi _{n},\quad x\in \partial \Omega.  \nonumber
\end{eqnarray}
Here, $A_{n}=(1+k_{n}/k_{n-1})$. We refer to \cite{KLN} for more details of
this derivation.

Recall that both functions $q$ and $V$ are unknown. Hence, keeping in mind
that we work with an iterative process, we have replaced $\nabla V,\Delta V$
with $\nabla V_{n-1},\Delta V_{n-1}$ in (\ref{12}), where $V_{n-1}$ is the
tail function calculated on the previous iteration. Also, we have replaced
the term $A_{n}\nabla q_{n}\nabla V_{n-1}$ with the term $A_{n}\nabla
q_{n-1}\nabla V_{n-1}$. In doing so, we have used $|\nabla q_{n}-\nabla
q_{n-1}|=O(h),\,h\rightarrow 0$ (see Lemma 4.3 of \cite{KLN}).

\subsection{The algorithm}

\label{sec:3.3}

In this algorithm, we use the index $n$ for the outer iterations. To update
the function $V$, we use, for each $n,$ inner iterations, which we number by
the index $i$.

\begin{center}
\textbf{Algorithm:}
\end{center}

\begin{itemize}
\item Set $q_{1,0}\equiv 0$, $\nabla V_{1,0}=\nabla V_{0}$, where $V_{0}$ is
the zero approximation for the tail function (subsection 3.4).

\item For $n=1$ to $N$,

\begin{enumerate}
\item Set $q_{n,0}=q_{n-1}$, $\nabla V_{n,0}=\nabla V_{n-1}$, $\Delta
V_{n,0}=\Delta V_{n-1}.$

\item For $i=1$ to $m$

\begin{itemize}
\item Find $q_{n,i}$ by solving the elliptic Dirichlet boundary value
problem (\ref{12}) where $V_{n-1}:=V_{n,i-1}$.

\item Compute $\nabla v_{n,i}=-(h\nabla q_{n,i}+h\nabla \overline{q_{n-1}}%
)+\nabla V_{n,i-1}$ and then $\Delta v_{n,i}=\func{div}\left( \nabla
v_{n,i}\right) $.

\item Compute $c_{n,i}$ \ using (\ref{eq:c_from_v}), 
\begin{equation}
c_{n,i}(x)=-\frac{1}{k_{n}^{2}}\left( \Delta v_{n,i}+(\nabla
v_{n,i})^{2}\right) .  \label{13}
\end{equation}

\item Find $u_{n,i}(x,\overline{k})$ by solving the Lippmann-Schwinger
equation (\ref{eq:lippmann-schwinger}) with the function $c\left( x\right)
:=c_{n,i}\left( x\right) -1$.

\item Update the gradient of the tail function, 
\begin{equation}
\nabla V_{n,i}(x)=\frac{\nabla u_{n,i}(x,\overline{k})}{u_{n,i}(x,\overline{k%
})}.  \label{eq:tail_grad}
\end{equation}
\end{itemize}

\item Set $q_{n}=q_{n,m}$, $c_{n}=c_{n,m}, V_n = V_{n, m}$.
\end{enumerate}

\item Let the pair of numbers $\left( n_{0},i_{0}\right) \in $ $[1,N]\times %
\left[ 1,m\right] $ be the one on which the stopping criterion is reached.
Set the function $c_{n_{0},i_{0}}(x)$ as the computed solution of the above
CIP, $c_{comp}\left( x\right) :=c_{n_{0},i_{0}}(x)$.
\end{itemize}

The stopping criterion for this iterative process should be developed
computationally, as it was also done in \cite{IPexp1,KLN,Liem,TBKF1}, see
Remark 3.1 at the end of subsection 3.4 for an explanation. We refer to
section 5.4 for some details regarding formula (\ref{13}).

\subsection{The first tail function}

\label{sec:3.4}

\label{sec:tail} The global convergence of the above algorithm is ensured by
the choice of the first tail function $V_{0}.$ It is important that we
choose the function $V_{0}$ without using any a priori knowledge of a small
neighborhood of the exact solution of our CIP.

Using an asymptotic expansion of the function $u\left( x,k\right) $ with
respect to $k$ as well as (\ref{eq:v}), it was proven in \cite{KLN} that,
under some conditions, there exists a function $p(x)\in C^{2}(\overline{%
\Omega })$ such that the following asymptotic behavior holds:


\begin{equation}
v(x,k)=ikp\left( x\right) \left( 1+O(1/k)\right) ,\quad k\rightarrow \infty
,x\in \overline{\Omega },  \label{14}
\end{equation}%
\begin{equation}
\left\vert O(1/k)\right\vert \leq B/k,\quad \forall x\in \overline{\Omega },
\label{15}
\end{equation}%
where $v(x,k)$ is defined in (\ref{eq:v}) and the number $B=B\left( c,\Omega
\right) >0$ depends only on the function $c$ and the domain $\Omega .$

By one of concepts of the regularization theory, one should assume that
there exists an idealized solution of the ill-posed problem under
consideration. This solution corresponds to the noiseless data, see, e.g. 
\cite{BK,T}. Let the function $c^{\ast }(x)$ be that exact solution. Below,
functions whose notations have the superscript \textquotedblleft $^{\ast }$"
denote those which are related to $c^{\ast }.$ Since the number $\overline{k}
$ is sufficiently large, then taking into account (\ref{15}), we drop the
term $O(1/k)$ in (\ref{14}) \ for $k\geq \overline{k}$. Hence, the function $%
v^{\ast }(x,k)$ is approximated as 
\begin{equation}
v^{\ast }(x,k)=ikp^{\ast }\left( x\right) ,\quad \forall k\geq \overline{k}.
\label{150}
\end{equation}%
Hence, we obtain 
\begin{equation}
V^{\ast }(x)=i\overline{k}p^{\ast }(x).  \label{16}
\end{equation}%
Since by (\ref{eq:q}) $q^{\ast }(x,k)=\partial _{k}v^{\ast }(x,k),$ then (%
\ref{150}) implies that 
\begin{equation}
q^{\ast }(x,k)=ip^{\ast }\left( x\right) ,\quad \forall k\geq \overline{k}.
\label{17}
\end{equation}%
Set in (\ref{9}) $k:=\overline{k}.$ Then substituting (\ref{16}) and (\ref%
{17}) in (\ref{9}) and using (\ref{10}) we obtain the Dirichlet boundary
value problem for the Laplace equation, 
\begin{eqnarray}
\Delta p^{\ast } &=&0,\quad x\in \Omega ,  \label{18} \\
p^{\ast } &=&i\psi ^{\ast }(x,\overline{k}),\quad x\in \partial \Omega .
\label{19}
\end{eqnarray}%
Hence, we set the first approximation $V_{0}(x)$ for the tail function as 
\begin{equation}
V_{0}(x)=i\overline{k}p(x),\quad x\in \Omega ,  \label{20}
\end{equation}%
where the function $p(x)$ is the solution of the following analog of the
problem (\ref{18}), (\ref{19})%
\begin{eqnarray}
\Delta p &=&0,\quad x\in \Omega ,  \label{180} \\
p &=&i\psi (x,\overline{k}),\quad x\in \partial \Omega .  \label{190}
\end{eqnarray}%
Assuming that functions $\psi ^{\ast }(x,\overline{k}),\psi (x,\overline{k}%
)\in C^{2+\alpha }\left( \partial \Omega \right) ,$ where $C^{2+\alpha
},\alpha \in \left( 0,1\right) $ is the H\"{o}lder space, it was proven in 
\cite{KLN} that 
\begin{equation}
\left\Vert V_{0}-V^{\ast }\right\Vert _{C^{2+\alpha }\left( \overline{\Omega 
}\right) }\leq C\overline{k}\Vert \psi -\psi ^{\ast }\Vert _{C^{2+\alpha
}(\partial \Omega )},  \label{1900}
\end{equation}%
where $C=C(\Omega )$ is a positive constant depending only on the domain $%
\Omega .$ Hence, the error in the first tail function depends only on the
error in the boundary data.

Therefore, we have obtained a good approximation of the tail function
already on the zero iteration of our algorithm. However, by our numerical
experience, we need to do more iterations. We point out that we use (\ref%
{150})-(\ref{190}) only on the zero iteration of our algorithm and do not
use them for other tail functions, which are calculated in our iterative
process. Also, equalities (\ref{150})-(\ref{19}) for the tail function
corresponding to the exact coefficient $c^{\ast }$ are used in the proof of
the global convergence theorem of \cite{KLN} only once: to obtain estimate (%
\ref{1900}), and they are not used for estimates of differences between
other iteratively obtained tail functions and the exact tail function $%
V^{\ast }.$

Equalities (\ref{150}) and (\ref{16}) represent our reasonable mathematical
assumption mentioned in Introduction. This assumption is justified by (\ref%
{14}) and (\ref{15}). We note that in (\ref{150}) and (\ref{16}) we do not
use any a priori information about a small neighborhood of the exact
coefficient $c^{\ast }.$ It was proven in \cite{KLN} that, given (\ref{150}%
)-(\ref{190}), the above algorithm delivers some points in a sufficiently
small neighborhood of the function $c^{\ast }(x).$ The approximation
accuracy depends only on the noise level in the data, the discretization
step $h$, the iteration number $n$ and the domain $\Omega .$ Therefore, our
algorithm converges globally within the framework of the assumption (\ref%
{150}), (\ref{16}).

\textbf{Remark 3.1. }We now explain why the stopping criterion for our above
algorithm should be chosen computationally. Although the global convergence
theorem of \cite{KLN} guarantees that this algorithm delivers some
coefficients, which are sufficiently close to the exact coefficient $c^{\ast
},$ it does not guarantee that the distances between the computed
coefficients and $c^{\ast }$ decrease with the iteration number. In other
words, even though those coefficients are located in that small
neighborhood, we do not know which of them is closer to $c^{\ast }$ than the
others. This explains why we should choose the stopping rule
computationally. We remind that the computed coefficient from our method can
always be refined by some locally convergent methods, see, e.g. \cite{BK}.
We also recall that the iteration number is often considered as the
regularization parameter in the theory of ill-posed problems \cite{BK,T}.

\section{Data acquisition and preprocessing}

\label{sec:4}

\subsection{Data acquisition}

\label{sec:4.1}

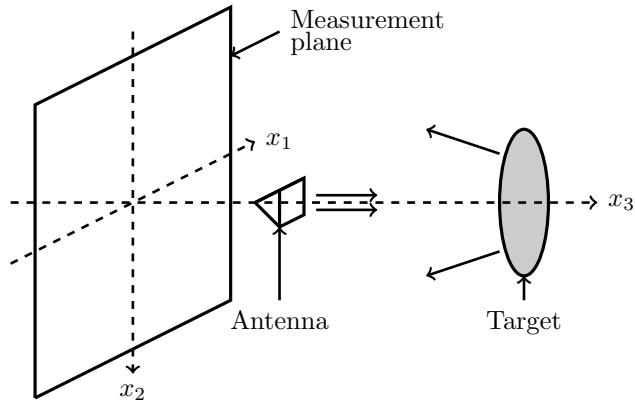
\begin{figure}[tbp]
\par
\begin{center}
\begin{tikzpicture}[font=\footnotesize,scale=0.65]
   
\draw [line width=1.2] (-2,-4) -- (2,-2) -- (2,4) -- (-2, 2) -- (-2,-4);
\draw[arrows=->,line width=1.0] (3,3.5) -- (2,3); \node [right] at (3,3.75) {Measurement}; \node [right] at (3,3.25) {plane}; 

\draw [line width=1.2] (3,-0.5) -- (3.5,-0.25) -- (3.5,0.5) -- (3,0.25) -- (3, -0.5);
\draw [line width=1.2] (2.5,0) -- (3,-0.5);
\draw [line width=1.2] (2.5,0) -- (3,0.25);
\draw [line width=1.2] (2.5,0) -- (3.5,0.5);
\draw[arrows=->,line width=1.0] (3,-2) -- (3,-0.5); \node [below] at (3,-2) {Antenna}; 

\draw [line width=1.2, fill=black!20!white]  (8,0) ellipse (0.5 and 1.5);
\draw[arrows=->,line width=1.0] (8,-2) -- (8,-1.5); \node [below] at (8,-2) {Target}; 

\draw[arrows=->,line width=1.0, dashed] (-2.5,-1.25) -- (2.5,1.25); \node [right] at (2.5,1.25) {$x_1$}; 
\draw[arrows=->,line width=1.0, dashed] (0,3.5) -- (0,-3.5); \node [below] at (0,-3.5) {$x_2$}; 
\draw[arrows=->,line width=1.0, dashed] (-2.5,0) -- (9.5,0); \node [right] at (9.5,0) {$x_3$}; 

\draw[arrows=->,line width=1.0] (3.75,-0.15) -- (5,-0.15); 
\draw[arrows=->,line width=1.0] (3.75,0.15) -- (5,0.15); 

\draw[arrows=->,line width=1.0] (7.5,1) -- (6,1.5);
\draw[arrows=->,line width=1.0] (7.5,-1) -- (6,-1.5);  

\end{tikzpicture}
\end{center}
\par
\caption{Schematic diagram of our experimental setup.}
\label{fig:geom}
\end{figure}

The experimental setup includes an emitting antenna, a detector, which is
moved on a vertical plane referred to as the measurement plane, and targets
to be imaged. The schematic diagram of our measurement arrangement is
presented on Figure \ref{fig:geom}. Consider the coordinate system $%
Ox_{1}x_{2}x_{3}$. Then $x_{1}$ and $x_{2}$ are horizontal and vertical axis
respectively. The $x_{3}-$axis is orthogonal to both $x_{1}$ and $x_{2}$
axis and points towards the target. Hence, the measurement plane, i.e. the
plane where detectors are located, is orthogonal to the $x_{3}-$axis. The
front face of each target is on the plane $\left\{ x_{3}=0\right\} .$
However, we do not use the information about the front face of the target
when applying our algorithm. The position of the emitting antenna is fixed,
i.e. antenna is not moving. Due to some technical issues, it was impossible
to place the antenna behind the measurement plane. Therefore, one among many
parasitic signals is due to the fact that the backscatter wave hits the
antenna and scatters again, also see Figure 3 for parasitic signals.

Below \textquotedblleft m" means meter. The distance between the measurement
plane and the emitting antenna is between 0.2 m and 0.25 m. Indeed, due to
some technical reasons, this distance was changing sometimes from one
experiment to another one. But we knew this distance in each experiment. The
distances between the measurement plane and the targets are about 0.8 m. As
it is shown in section 4.2.7, the minimal optimal frequency of our signal is
6.56 GHz (it varies from one target to another one). Hence, we have that the
maximal optimal wavelength $\lambda _{\max }=0.0457$ m. Hence, the
antenna/target distance of about $0.6$ m approximately equals to $%
13.1\lambda _{\max }$, which can be considered as a large distance in
Physics. \ Therefore, on such large distances we can approximate the
spherical wave emitted by the antenna as the plane wave. This justifies our
modeling of the incident signal by the plane wave, see (\ref{3}).

The detector is moved over a 1\ m $\times $ 1\ m square with the step size
of 0.02 m, starting from the top-left corner $(x_{1},x_{2})=(-0.5,-0.5)$ m
and ending at the bottom-right corner $(x_{1},x_{2})=(0.5,0.5)$ m. Again,
this square is referred to as the measurement plane. At each position of the
detector, the antenna emits a number of electric pulses with duration of 300
picoseconds (ps) each and the detector records signals (voltage) during 10
nanoseconds time period with the time step $\Delta t=10$ picoseconds (ps). A
part of the emitted signal goes forward towards the target and another part
goes backwards. We call the second part \textquotedblleft direct signal".
Recorded signals contain both the direct signals from the antenna and the
backscatter signals both from the targets and surrounding objects. To reduce
instabilities, each pulse is repeated 800 times at each position of the
detector and the measured data on this detector is averaged then over these
800 signals.

Now we present our setup in the dimensionless variables. Denote by $%
c_{0}\approx 0.0003$ (m/ps) the speed of light in air. In the case of
dimensions, equation (\ref{4}) is given by 
\begin{equation}
\frac{c\left( x\right) }{c_{0}^{2}}w_{tt}=\Delta w.  \label{time}
\end{equation}%
%
%
%
%
%
We introduce dimensionless variables for the latter equation as 
\[
y=\frac{x}{0.1\mathrm{m}},\quad t^{\prime }=\frac{0.0003}{\mathrm{0.1ps}}t. 
\]%
Then substitution of $c_{0}$ in (\ref{time}) and a straightforward
calculation lead to %
\[
c\left( y\right) w_{t^{\prime }t^{\prime }}(y,t^{\prime })=\Delta
_{y}w(y,t^{\prime }), 
\]%
which is the dimensionless equation (\ref{4}). We now can apply the Fourier
transform 
\[
\widetilde{u}\left( y,k\right) =\int_{0}^{\infty }w\left( y,t\right)
e^{ikt}dt 
\]%
with dimensionless $k$ like in (\ref{eq:helmholtz}). It follows from the
results of \cite{V} that the function $\widetilde{u}$ satisfies conditions (%
\ref{eq:helmholtz})-(\ref{eq:radiation}). In particular, 
\[
\Delta _{y}\widetilde{u}+k^{2}c\left( y\right) \widetilde{u}=0. 
\]%
For brevity, we keep the same notations for dimensionless variables. From
now on, \textquotedblleft 5" of the dimensionless length means 0.5 m of
length. Now from a simple calculation, the relation between the
dimensionless wave number $k$ and the frequency $f$ (in Hertz) is given by $%
f=c_{0}k/(2\pi \times 0.1m)$. %
Therefore, the optimal dimensionless wave number corresponding to the
frequency $6.56$ GHz is 
\begin{equation}
k_{opt}=13.75.  \label{21}
\end{equation}


\subsection{Data preprocessing}

\label{sec:4.2}

As it was pointed out in Introduction, our experimental data are
substantially different from the computationally simulated data.\ 

\subsubsection{Main indicators of the huge discrepancy between experimental
and computationally simulated data}

\label{sec:4.2.1}

These indicators are:

\begin{enumerate}
\item[1.] The backscatter signal is unstable due to the instability of the
emitted signal.

\item[2.] The emitted signal moves not only forwards towards the target but
also moves backwards to the measurement plane. Recall that we call the
latter the \textquotedblleft direct signal".

\item[3.] There were some metallic parts of our device placed behind the
measurement plane. So, both direct and backscatter signals are scattered by
these parts and come back to the detector.

\item[4.] The emitting antenna is placed between the measurement plane and
the targets. This causes the shadow effect on the measurement plane as well
as some additional parasitic signals.

\item[5.] Our signals are scattered by some objects located in the room,
such as furniture, and these scattering signals are also recorded by the
detector.

\item[6.] Scales of magnitudes of experimental and computationally simulated
data are significantly different.

\item[7.] The Fourier transform of the preprocessed time dependent data has
the so-called \textquotedblleft optimal frequency". The data are reliable
only in a small interval of frequencies surrounding the optimal one and are
reliable on the rest of the frequency interval. The same was observed for
the case of multi-frequency experimental data in \cite{Liem}. However, this
phenomenon is not the case in the computationally simulated data.
\end{enumerate}

The above leads to the conclusion that it is necessary to apply a heuristic
data preprocessing procedure. This procedure aims to make the processed data
look at least somehow similar to the computationally simulated data. The
preprocessed data are then used as the input for the above algorithm, i.e.
as the function $g\left( x,k\right) $ in (\ref{eq:cip_fdd}). Below in this
subsection 4.2 we describe steps of our data preprocessing procedure. This
procedure consists of two stages. First, we preprocess the data in the time
domain. Second, we apply the Fourier transform to the preprocessed data and
then we preprocess those Fourier transformed data. When working in the time
domain, our procedure is similar to the one of \cite{TBKF1} and when
preprocessing the data in the frequency domain, our procedure is similar to
the one of \cite{Liem}.

However, there are some important differences from the preprocessing
procedures of \cite{Liem,TBKF1}. On the first stage the difference with \cite%
{TBKF1} is in Step 3. We extract signals from the target relying on the
estimates of travel time of the measured signals. This method turned out to
be simpler than the one used in \cite{TBKF1}. 
On the second stage, the most important difference with \cite{Liem} is that,
unlike \cite{Liem}, we \textquotedblleft shift" the data with respect to the
frequency, see (\ref{22}) below. We need to do this shift in order to
decrease our computational burden. Such a shift was not necessary to do in 
\cite{Liem}, see section 4.2.7. Steps 1-3 below are of stage 1 and steps 4-7
are of stage 2.

\subsubsection{Step 1. Off-set correction}

\label{sec:4.2.2}

The measured signals are usually shifted from the zero mean value, so we
subtract the mean value from them.

\subsubsection{Step 2. Time-zero correction}

\label{sec:4.2.3}

The moment, when the incident pulses are emitted from the antenna, is
referred to as the time-zero. In this step, using the direct signal from the
antenna, the data is shifted to the correct time-zero. Time-zero correction
is essential for extracting useful signals.

\subsubsection{Step 3. Extraction of scattered signals}

\label{sec:4.2.4}

As mentioned above, our measured data contain various types of parasitic
signals. However, we need to use only the signals scattered by the targets.
In our experimental setup we know the distance between the measurement plane
and the antenna. In addition, we can roughly estimate the distances between
the measurement plane and targets. These lead to estimates of the travel
time of the signal between the antenna and the targets, which, in turn
results in estimates of the travel time between the targets and each
location of the detector on the measurement plane. Next, all signals before
this moment of arrival are eliminated. Note that the travel time is
different for each location of the detector on the measurement plane. In
addition, we know an estimate from the above of the linear sizes of any
target. In fact, these sizes should approximately match sizes of
antipersonnel mines and IEDs. So, each linear size of each target is between
0.05 m and 0.15 m. Thus, this estimate helps us to eliminate useless signals
which are recorded after the true signals scattered by targets are recorded.

\begin{figure}[tbp]
\centering
\subfloat[\label{fig:dataset_a}]{\includegraphics[width=0.49%
\textwidth]{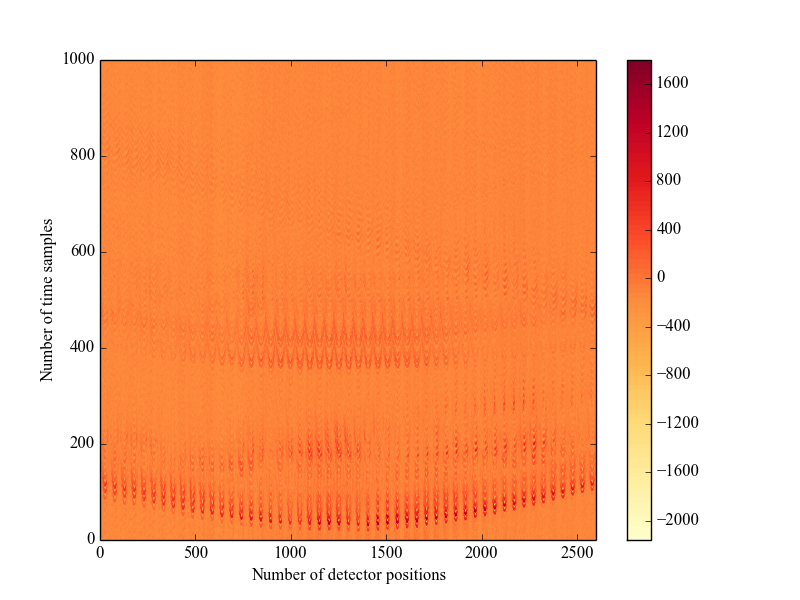}} \subfloat[\label{fig:dataset_b}]{%
\includegraphics[width=0.49\textwidth]{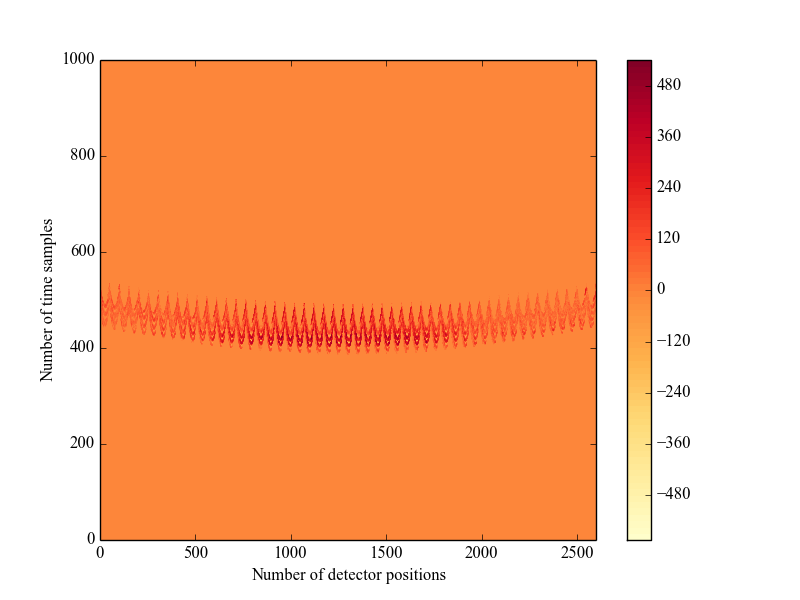}}
\caption{Experimental data before (a) and after steps 1-3 of the data
prepocessing procedure (b).}
\label{fig:dataset}
\end{figure}

The experimental data for one of our targets recorded before and after the
previous three steps of preprocessing are depicted on Figure \ref%
{fig:dataset}. Here, the data is depicted in the form of a matrix, whose
rows represent positions of the detector and columns represent time samples.
We see in Figure \ref{fig:dataset}-a) that parasitic signals in the raw data
dominate useful ones. Next, Figure \ref{fig:dataset}-b) shows that parasitic
signals were successfully excluded after the above step 3 of data
preprocessing.

\begin{figure}[tbp]
\par
\begin{center}
\includegraphics[width=0.7\textwidth]{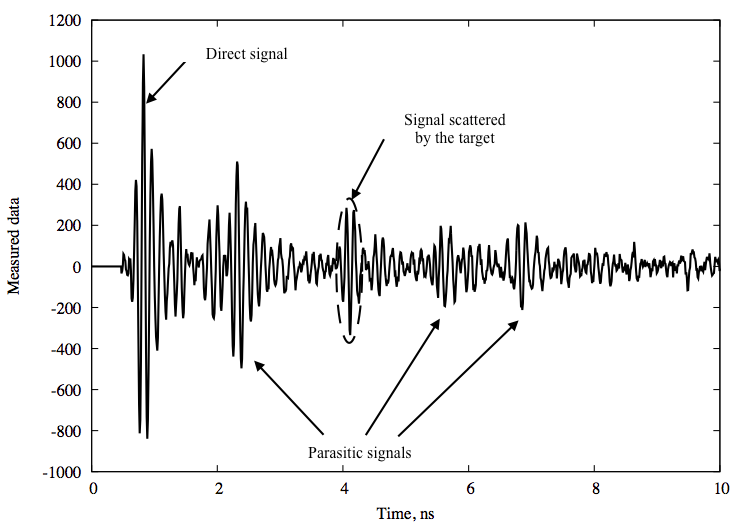}
\end{center}
\caption{A typical sample of our time dependent experimental data.}
\label{fig:time_curve}
\end{figure}

We now show how to actually perform Step 3. As an example, we plot in Figure %
\ref{fig:time_curve} the time dependent data for a target on a detector
located in the middle of the measurement plane. We know that the first
strong signal in this curve is the direct signal from the antenna. We also
know that the second strong signal after the direct one is from objects
located behind the measurement plane. Our estimate of the distance between
the front face of the target and the detector enables us to estimate the
time of arrival of the signal scattered by the front face the target, i.e.
the travel time. Hence, we set to zero the part of the curve which is before
that time of arrival. Next, our estimate of the linear size of the target
allows us to estimate the time of arrival of the signal, which is scattered
by the back side of the target. The part of the curve, which is after the
latter signal, is set to zero. Figure 4 displays the preprocessed curve of
Figure 3 after Steps 1-3.

\begin{figure}[tbp]
\par
\begin{center}
\includegraphics[width=0.7\textwidth]{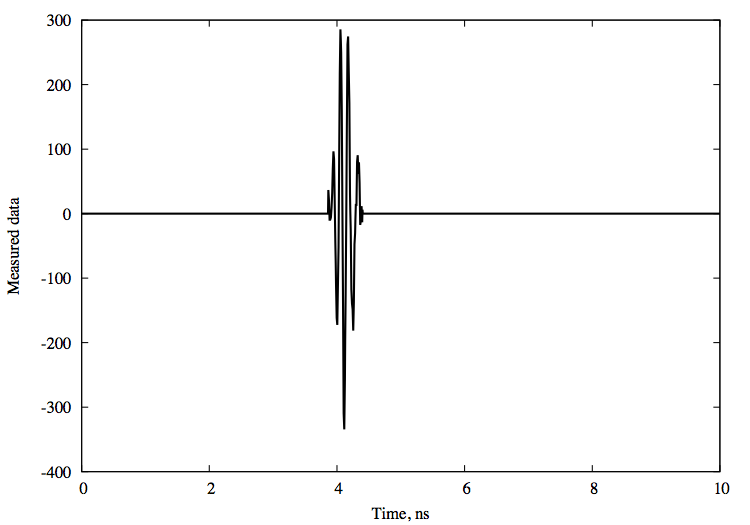}
\end{center}
\caption{Preprocessed signal of Figure 3 after Steps 1-3.}
\label{fig:time_curve2}
\end{figure}

\subsubsection{Step 4. Fourier transform}

\label{sec:4.2.5}

Let $f(x,t)$ be the time dependent data obtained after above steps 1-3 of
the data preprocessing. Following (\ref{7}), we calculate the Fourier
transform with respect to $t$, 
\[
\varphi (x,k)=\int_{0}^{\infty }f(x,t)\,e^{ikt}\,dt,\quad k\in \lbrack 
\underline{k},\overline{k}]. 
\]%
Since the data are measured on the rectangle $S_{m}=\left( -5,5\right)
\times \left( -5,5\right) $ in the measurement plane, denoted by $P_{m}$, we
set 
\begin{equation}
f(x,t)=\varphi (x,k)=0,\quad \left( x_{1},x_{2}\right) \not\in S_{m}
\label{100}
\end{equation}%
Figure \ref{fig:f} shows $|\varphi (x,k)|$ $(k=13.5,$ see section 4.2.7$)$
for the experimental data for one of our targets on the measurement plane.
Note that since the function $f(x,t)$ is extracted from the backscatter
data, then the function $\varphi (x,k)$ also corresponds to the backscatter
data rather than to a full wave field.

\subsubsection{Step 5. Data propagation}

\label{sec:4.2.6}

The measurement plane was located quite far from the targets. This means
that we were supposed to solve the CIP in a large domain, which is
inconvenient. Besides, stability estimates for solving boundary value
problems (\ref{12}) and (\ref{180}), (\ref{190}) are worsened for large
domains, as compared with smaller ones. The goal of the data propagation
procedure is twofold: (a) To reduce the computational domain and (b) To make
the data more \textquotedblleft clear", so that they the $x_{1},x_{2}$
coordinates of the targets would be clarified. The data propagation
procedure is described in details in \cite{Nov,TBKF2}. Nevertheless, for the
convenience of the reader, we briefly present the main steps below. This
procedure is an application of the angular spectrum representation which is
a well-known technique in the optics community, see \cite{Born, Nov}.
Basically, this procedure \textquotedblleft moves" the data closer to the
targets. The idea is that we consider the Helmholtz equation for the
backscatter wave field outside the scattering medium. Then the Fourier
transform with respect to $x_{1}$ and $x_{2}$ of the solution of this
equation satisfies an ordinary differential equation in the $x_{3}$%
-direction. Together with the radiation condition and the boundary condition
on the measurement plane one can solve that 1D problem and obtain the data
on the propagated plane.

Let the measurement plane be $P_{m}=\{x_{3}=b\}$, $b>0$ and let $%
P_{p}=\{x_{3}=a\},a<b$ be the plane to which we want to propagate our data.
We need to determine the function $u_{sc}\left( x,k\right) =\tilde{g}(x,k)$
on $P_{p},$ given the data $\varphi (x,k)$ on $P_{m}$. The computation
consists of two steps. First, we compute the Fourier transform with respect
to $x_{1}$ and $x_{2}$ of the function $\varphi (x,k)$ for $x\in P_{m},$ 
\[
p(k_{1},k_{2},k)=\frac{1}{2\pi }\displaystyle\int\limits_{\mathbb{R}%
^{2}}\varphi (x,k)\,\exp \left( i\left( x_{1}k_{1}+x_{2}k_{2}\right) \right)
\,dx_{1}dx_{2}, 
\]%
also see (\ref{100}). In the second step, we apply the inverse Fourier
transform, where the integration is carried out over a bounded domain $%
\left\{ k_{1}^{2}+k_{2}^{2}<k^{2}\right\} .$ Denote $k_{z}=\sqrt{%
k^{2}-k_{1}^{2}-k_{2}^{2}}>0$. Then 
\[
u_{sc}\left( x,k\right) =\tilde{g}(x,k)=\frac{1}{2\pi }\int%
\limits_{k_{1}^{2}+k_{2}^{2}<k^{2}}p(k_{1},k_{2},k)%
\,e^{-i(x_{1}k_{1}+x_{2}k_{2}+(b-a)k_{z})}dk_{1}dk_{2},x_{3}=a. 
\]

We propagate our experimental data for all targets from the measurement
plane $P_{m}=\{x_{3}=-8\}$ to the propagated plane $P_{p}=\{x_{3}=-0.75\}$.
The functions $\left\vert \varphi (x,k)\right\vert $ and $\left\vert \tilde{g%
}(x,k)\right\vert $ for both simulated and experimental data at $k=13.5$
(see section 4.2.7) are displayed on Figures 5. This is the case of a single
target, which is the target number 1, see Table 1 below. Figure 6 displays
functions $\left\vert \varphi (x,k)\right\vert $ and $\left\vert \tilde{g}%
(x,k)\right\vert $ for the case of two inclusions for $k=15.1$ (see section
4.2.7). Most importantly, the data on the propagated plane are focused near
the location of the target rather than being spread out as on the
measurement plane. One can see on Figures 5b and 6b that we can now estimate 
$x_{1},x_{2}$ coordinates of targets, unlike Figures 5a and 6a. Furthermore,
unlike the measurement plane of Figure 6a, two inclusions are clearly
separated on the propagated plane of Figure 6b. Using these observations, we
set%
\begin{equation}
\tilde{g}(x,k)=0,\quad \mathrm{\ for}\ \left( x_{1},x_{2}\right)
\notin S_{p}=\left( -2.5,2.5\right) \times \left( -2.5,2.5\right) \subset
P_{p}.  \label{102}
\end{equation}

\subsubsection{Step 6. Choosing the interval for wave numbers}

\label{sec:4.2.7}

We now use temporary notations for numbers $\underline{k}$ and $\overline{k}$
as $\underline{k}^{\prime }$ and $\overline{k}^{\prime }$ respectively. We
explain in this section how do we choose an appropriate interval $[%
\underline{k}^{\prime },\overline{k}^{\prime }]$ of wave numbers for each
target. The interval of wave numbers should meet the following criteria:

\begin{enumerate}
\item The interval $[\underline{k}^{\prime },\overline{k}^{\prime }]$
contains the global maximizer of the function $s(k),$ 
\begin{equation}
s(k)=\max_{x\in S_{p}}|\tilde{g}(x,k)|.  \label{1010}
\end{equation}%
Indeed, it seems to be intuitively at least that the data are most reliable
for those values of wave numbers which are near that maximizer.

\item For $k\in \lbrack \underline{k}^{\prime },\overline{k}^{\prime }]$,
the propagated data indicate correctly $x_{1},x_{2}$ coordinates of targets.
\end{enumerate}

Here is an example. We plot on Figure 7a the graph of the function $s(k)$
defined in (\ref{1010}) for the target number 1 (see Table \ref{tab:targers}%
). This is our reference target for data calibration (see section 4.2.8).
Since this is the reference target, we are supposed to know its location,
shape and the value of its dielectric constant. We point out that this is
unlike all other targets.

We observe on Figure 7a that $\max_{\left[ 0,20\right] }s\left( k\right)
=s\left( 13.75\right) .$ Hence, the wave number $k=13.75$ is referred as the
optimal frequency for the target number 1, see (\ref{21}). Thus, we choose
the interval of wave numbers for the target number 1 as $[\underline{k}%
^{\prime },\overline{k}^{\prime }]=[13.5,14]$ so that the optimal frequency
is in the middle of this interval, $k_{opt}^{\prime }=\left( \underline{k}%
^{\prime }+\overline{k}^{\prime }\right) /2$. We choose the length of the $%
k- $interval to be the same for all targets, i.e. $a=\overline{k}^{\prime }-%
\underline{k}^{\prime }=0.5$. Figure 7b also confirms that this interval of
wave numbers is optimal. It is clear from Figure 7b that the absolute values
of the propagated data for $k\in \lbrack 13.5,14]$ are near the maximal
value of the function $|\tilde{g}(x,k)|$ and these data strongly focus in a
subdomain of the rectangle $S_{p}$.

The interval $[13.5,14]$ of wave numbers corresponds to the interval of
frequencies $\left[ 6.51,6.69\right] $ GHz. This can be converted to the
interval of wavelengths $\lambda \in \left[ 0.0448,0.046\right] $ m. Note
that the wave number interval for each target is different. The optimal
frequency changes from 13.75 for the target number 1 to 18.7 for target
number 9,%
\begin{equation}
k_{opt}\in \left[ 13.75,18.7\right] .  \label{23}
\end{equation}



It is important to observe on Figure 7a that the function $s\left( k\right) $
changes slowly in the interval $[\underline{k}^{\prime },\overline{k}%
^{\prime }]=[13.5,14]$. On the other hand, it changes very rapidly for $k\in
\left( 11,\underline{k}^{\prime }\right) \cup \left( \overline{k}^{\prime
},16\right) ,$ which is outside of the interval $[\underline{k}^{\prime },%
\overline{k}^{\prime }]$.\ Indeed, it changes between 200 and 1800 for $k\in
\left( 11,\underline{k}^{\prime }=13.5\right) $ and it changes between 1750
and 600 for $k\in \left( \overline{k}^{\prime }=14,16\right) $. Hence, we
estimate the modulus of the derivative $\left\vert s^{\prime }(k)\right\vert 
$ as $\left\vert s^{\prime }(k)\right\vert >500$ for $k\in \left( 11,%
\underline{k}^{\prime }\right) \cup \left( \overline{k}^{\prime },16\right)
. $ We consider this as an indication that the data are possibly highly
unstable for $k\in \left( 11,\underline{k}^{\prime }\right) \cup \left( 
\overline{k}^{\prime },16\right) .$ Similar observations are true for the
interval $k\in \left[ 16,20\right] ,$ which surrounds the local maximizer of 
$s\left( k\right) $ at $k=18.$ However, since $s\left( 18\right) <s\left(
13.75\right) ,$ then the first above criterion is not satisfied for $k\in %
\left[ 16,20\right] .$ As to $k\in \left( 0,11\right) ,$ we have observed
that the second above criterion is not satisfied for our reference target
number 1. Hence, we have ignored the interval $k\in \left( 0,11\right) $ for
this target. Also, because of this, we did the same for similar intervals
for other targets without even checking the second criterion for them since
we are not supposed to know locations of targets other than of the reference
target. Thus, we work only on the interval $k\in \lbrack \underline{k}%
^{\prime },\overline{k}^{\prime }]=[13.5,14]$ for the target number 1 and on
similar intervals for other targets.

However, the selected interval $[\underline{k}^{\prime },\overline{k}%
^{\prime }]=[13.5,14]$ of wave numbers is inconvenient for our computational
purpose. The reason of this is that the solution $u\left( x,k\right) $ of
the forward problem (\ref{eq:helmholtz})-(\ref{eq:radiation}) is highly
oscillatory for $k\in \lbrack 13.5,14],$ so as for other intervals selected
for other targets, see (\ref{23}). This means that our solver for the
Lippmann-Schwinger equation (\ref{eq:lippmann-schwinger}) requires a very
fine mesh for these values of $k$ and thus, this solver becomes very
expensive in terms of both the computing time and memory. On the other hand,
we need to solve equation (\ref{eq:lippmann-schwinger}) on each iterative
step of our algorithm to update the function $u\left( x,\overline{k}\right)
. $

Hence, we act as follows: In all our computations we use $k\in \lbrack 
\underline{k},\overline{k}]=[6,6.5]$ for all targets. However, for the data
for the function $u_{sc}\left( x,k\right) $\ on the propagated plane $P_{p}$
we use%
\begin{equation}
u_{sc}\left( x,k\right) :=\tilde{g}(x,k^{\prime }),k^{\prime }=k^{\prime
}\left( k\right) =k+\underline{k}^{\prime }-6\in \lbrack \underline{k}%
^{\prime },\overline{k}^{\prime }],x\in S_{p},  \label{22}
\end{equation}%
also see (\ref{102}). Thus, (\ref{22}) means that we shift
the data from the interval $[\underline{k}^{\prime },\overline{k}^{\prime }]=%
\left[ \underline{k}^{\prime },\underline{k}^{\prime }+0.5\right] $ to the
interval $[\underline{k},\overline{k}]=[6,6.5].$ In other words, even though
the true data we work with are on the interval $[\underline{k}^{\prime },%
\overline{k}^{\prime }]$ of wave numbers, we pretend that they are on the
interval $k\in \lbrack 6,6.5].$

We have selected this $k-$interval since the same interval $[\underline{k},%
\overline{k}]=[6,6.5]$ was successfully used in the paper \cite{Liem} of our
group. The data shift was not performed in \cite{Liem}, since the optimal
wave number of the data of \cite{Liem} was inside the interval $[\underline{k%
},\overline{k}]=[6.25,6.70]$ for all targets. The latter corresponds to the
interval of frequencies $\left[ 2.98,3.19\right] $ GHz.

\textbf{Remark 4.1}. Following (\ref{11}), we have subdivided the interval $%
[6,6.5]$ in ten (10) subintervals with the step size $h=0.05$ and $N=10$. We
point out that, in the case of the target number 1, even though the function
(\ref{1010}) changes slowly on the interval $k\in \lbrack 13.5,14],$ the
function $|\tilde{g}(x,k)|$ might change rapidly with respect to $k$ for
this $k-$interval at some points $x$. But at least the above considerations
indicate that it seems to be meaningless to work on the intervals $k\in
\left( 0,13.5\right) $ and $k\in \left( 14,16\right) .$ The same procedure
of the choice of the interval of wave numbers and of the data shift was
performed for all targets used in our experiments. It follows from (\ref{23}%
) that optimal intervals for different targets were different. See section 7
for our ultimate judgement of the validity of the entire data preprocessing
procedure.

\subsubsection{Step 7. Data calibration}

\label{sec:4.2.9}

\begin{figure}[tbp]
\centering
\subfloat[\label{fig:f}]{\includegraphics[width=0.49\textwidth]{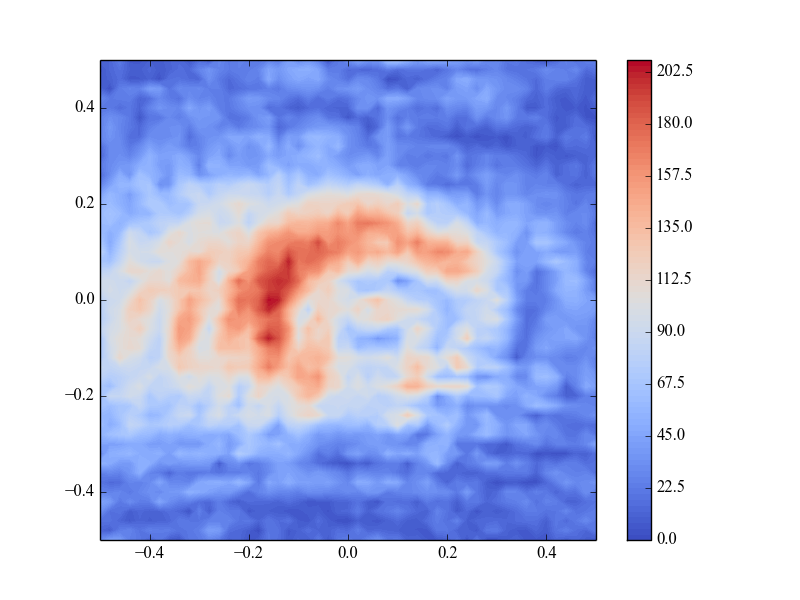}} %
\subfloat[\label{fig:g}]{\includegraphics[width=0.49\textwidth]{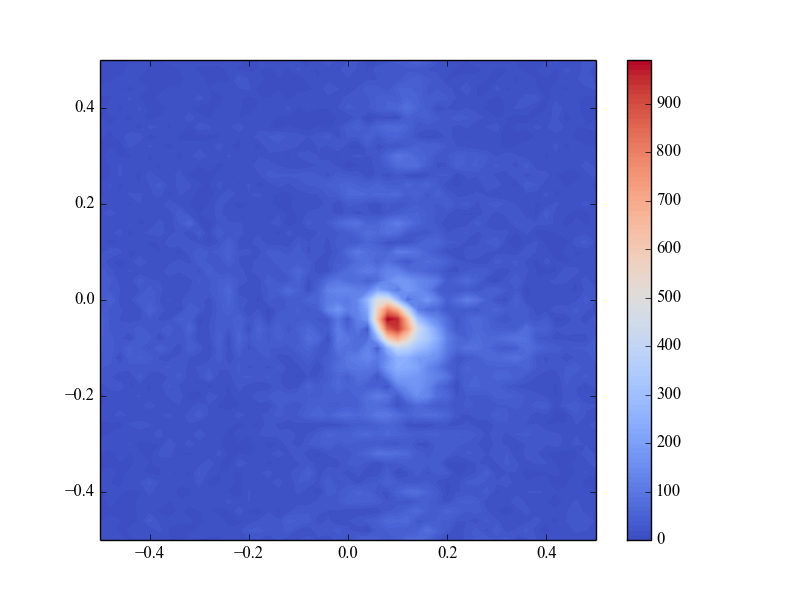}} 
\newline
\subfloat[\label{fig:f_sim}]{\includegraphics[width=0.49%
\textwidth]{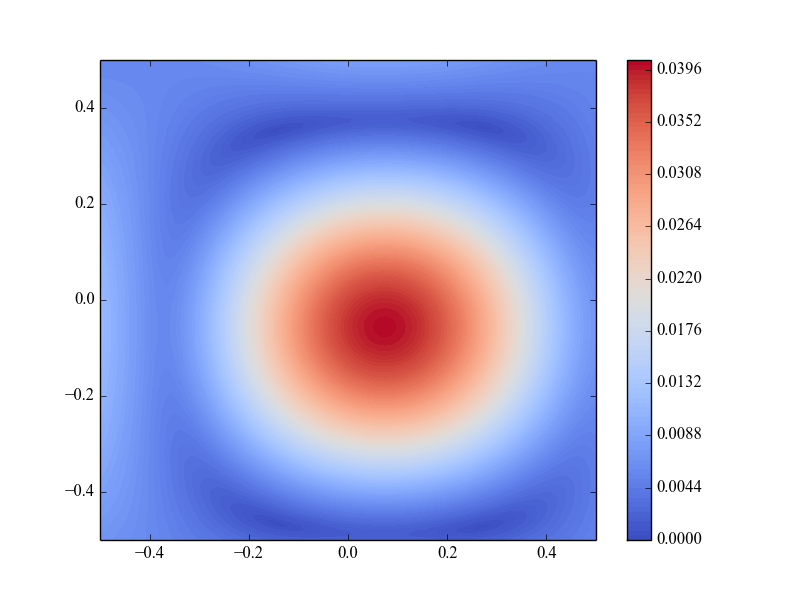}} \subfloat[\label{fig:g_sim}]{%
\includegraphics[width=0.49\textwidth]{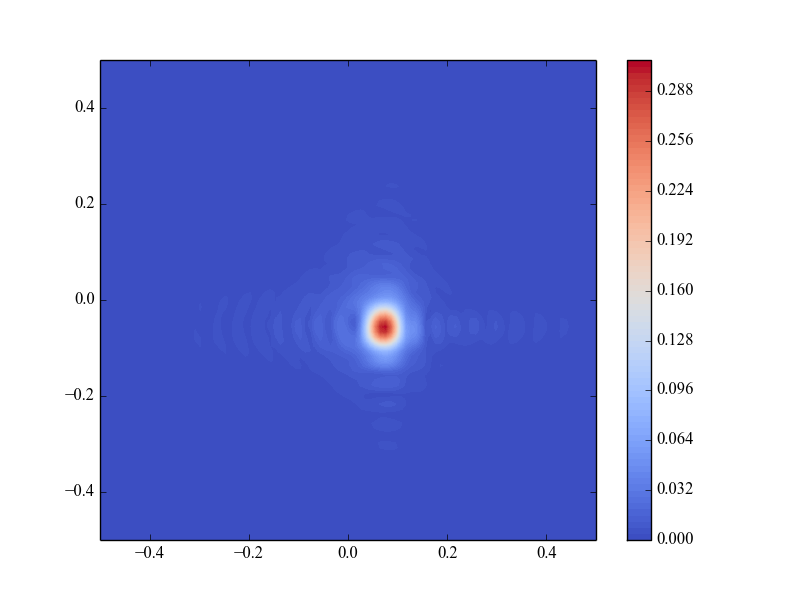}}
\caption{The absolute value of the Fourier transform of the preprocessed
time dependent experimental data on the measurement plane (a) and the
propagated plane (b). The absolute value of the Fourier transform of the
computationally simulated data for the same target on the measurement plane
(c) and on the propagated plane (d). $x_{1},x_{2}$ coordinates of the
inclusion are clearly seen on (b) and (d), unlike
(a) and (c). Here, $k= 13.5$.}
\label{fig:fourier_propag}
\end{figure}

We have observed that magnitudes of the experimental data are significantly
different from magnitudes of computationally simulated data. For example,
compare the absolute value of the Fourier transform of the experimental data
on the propagated plane (Figure \ref{fig:g}) with the absolute value of
propagated computationally simulated data on that plane (Figure \ref%
{fig:g_sim}). The maximal value of propagated experimental data is about
950, whereas the maximal value of the propagated simulated data is about
0.28. Therefore, we need to scale the experimental data to the same scale as
the one in computationally simulated data.

To do this, we consider one target, which we call \textquotedblleft
calibration target" (the target number 1) and consider the function $\tilde{g%
}(x,k+\underline{k}^{\prime }-6)$ for this target. Next, we compute the
function $\varphi \left( x,k\right) :=u\left( x,k\right) \mid _{x\in
S_{m}\subset P_{m}},k\in \lbrack \underline{k},\overline{k}]=[6,6.5]$ (also,
see (\ref{100})) on the measurement plane $P_{m}$ for a target, which mimics
the target number 1: it has the same size, location and dielectric constant
as the one of the target number 1. We compute the function $\varphi \left(
x,k\right) $ via the numerical solution of the Lippmann-Schwinger equation (%
\ref{eq:lippmann-schwinger}) for this target. Next, we propagate the
simulated data $\varphi \left( x,k\right) $ to the propagation plane $P_{p}$
and obtain the function $\tilde{g}_{sim}\left( x,k\right) $ this way.

Next, we define the $k-$dependent calibration factor $d(k)$ as 
\begin{equation}
d(k)=\frac{\max_{S_{p}}(|\tilde{g}_{sim}(x,k)|)}{\max_{S_{p}}(|\tilde{g}(x,k+%
\underline{k}^{\prime }-6)|)},k\in \lbrack \underline{k},\overline{k}%
]=[6,6.5].  \label{eq:calib_factor}
\end{equation}%
see (\ref{102}) and (\ref{22}). Once chosen, this function $d(k)$ remains
the same for other targets. Using (\ref{eq:calib_factor}), we set for all
targets%
\begin{eqnarray}
u\left( x,k\right) &:&=g\left( x,k\right) =\exp \left( ikx_{3}\right)
+d\left( k\right) \tilde{g}(x,k+\underline{k}^{\prime }-6),  \label{1000} \\
x &\in &S_{p}\subset P_{p},k\in \lbrack \underline{k},\overline{k}]=[6,6.5].
\nonumber
\end{eqnarray}%
In this formula we multiply the experimental backscatter data by the
calibration factor and then sum up with the incident wave. Note that it is
not necessary of course that the calibration factor $d(k)$ defined in (\ref%
{eq:calib_factor}) using only one target would be perfect for other targets.
Nevertheless, our reconstruction results below indicate that it works well.
Similar considerations about the calibration factor are true for other
studies of experimental data via GCM \cite{IPexp1,Liem,TBKF1,TBKF2}.

\begin{figure}[tbp]
\centering
\subfloat[\label{fig:f2}]{\includegraphics[width=0.49%
\textwidth]{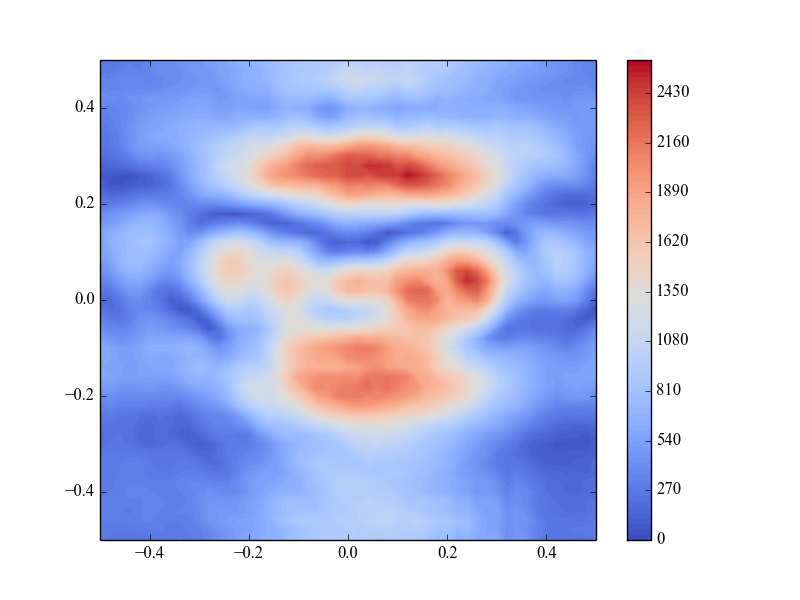}} \subfloat[\label{fig:g2}]{%
\includegraphics[width=0.49\textwidth]{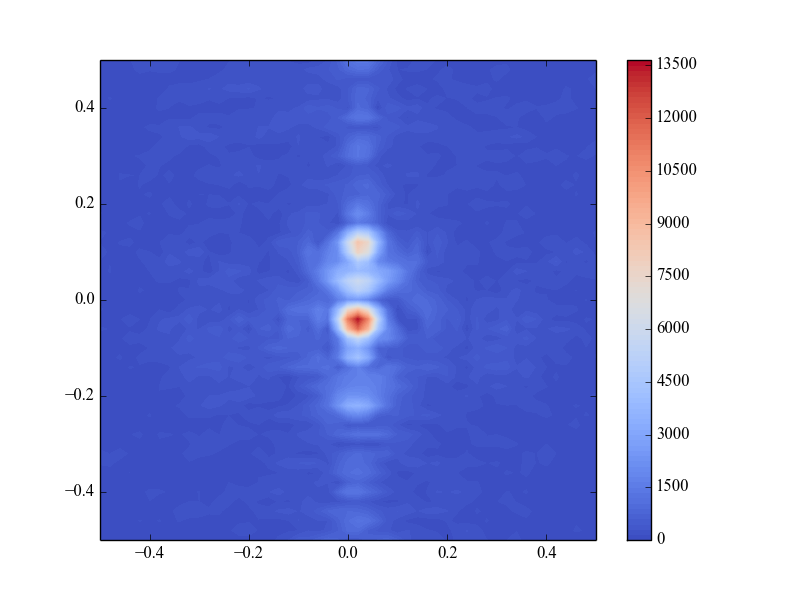}}
\caption{The absolute value of the Fourier transform of the experimental
data for the case of two inclusions on the measurement plane (a) and on the
propagated plane (b). $x_{1},x_{2}$ coordinates of both inclusions are
clearly seen on (b), unlike (a). Here, $k= 15.1$.}
\label{fig:fourier_propag_2inc}
\end{figure}

\section{Numerical implementation}

\label{sec:5}

In this section, we describe some details of the numerical implementation of
the algorithm described in section 3.2.

\subsection{Computational domain}

\label{sec:5.1}

Recall that the $x_{3}$-axis in our coordinate system points from the
measurement plane towards the target and the front face of the target is on
the plane $\left\{ x_{3}=0\right\} .$ We also recall that our propagated
plane is $P_{p}=\left\{ x_{3}=-0.75\right\} .$ We choose the computational
domain $\Omega $ as: 
\begin{equation}
\Omega =\{x\in (-2.5,2.5)\times (-2.5,2.5)\times (-0.75,4.25)\}.  \label{240}
\end{equation}%
The choice (\ref{240}) means that the domain $\Omega $ is a cube whose side
is 0.5 m. This is sufficient for our desired application since the linear
sizes of antipersonnel mines and IEDs usually are between about 5 and 15
centimeters (cm). Observe that by the choice (\ref{240}) we restrict our
attention in planes, which are orthogonal to $x_{3},$ to rectangles $\left(
x_{1},x_{2}\right) \in (-2.5,2.5)\times (-2.5,2.5),$ whereas the measured
data are given on the rectangle $\left( x_{1},x_{2}\right) \in
S_{m}=(-5,5)\times (-5,5)\subset P_{m}.$ This is because our propagated data
are set to zero outside of the rectangle $S_{p}\subset P_{p},$ see (\ref{102}%
). We solve Dirichlet boundary value problems (\ref{12}) by the finite
element method using the software FreeFem++ \cite{Hecht2012}.

\begin{figure}[tbp]
\par
\begin{center}
\subfloat[\label{fig 7b}]{\includegraphics[width=0.4\textwidth]{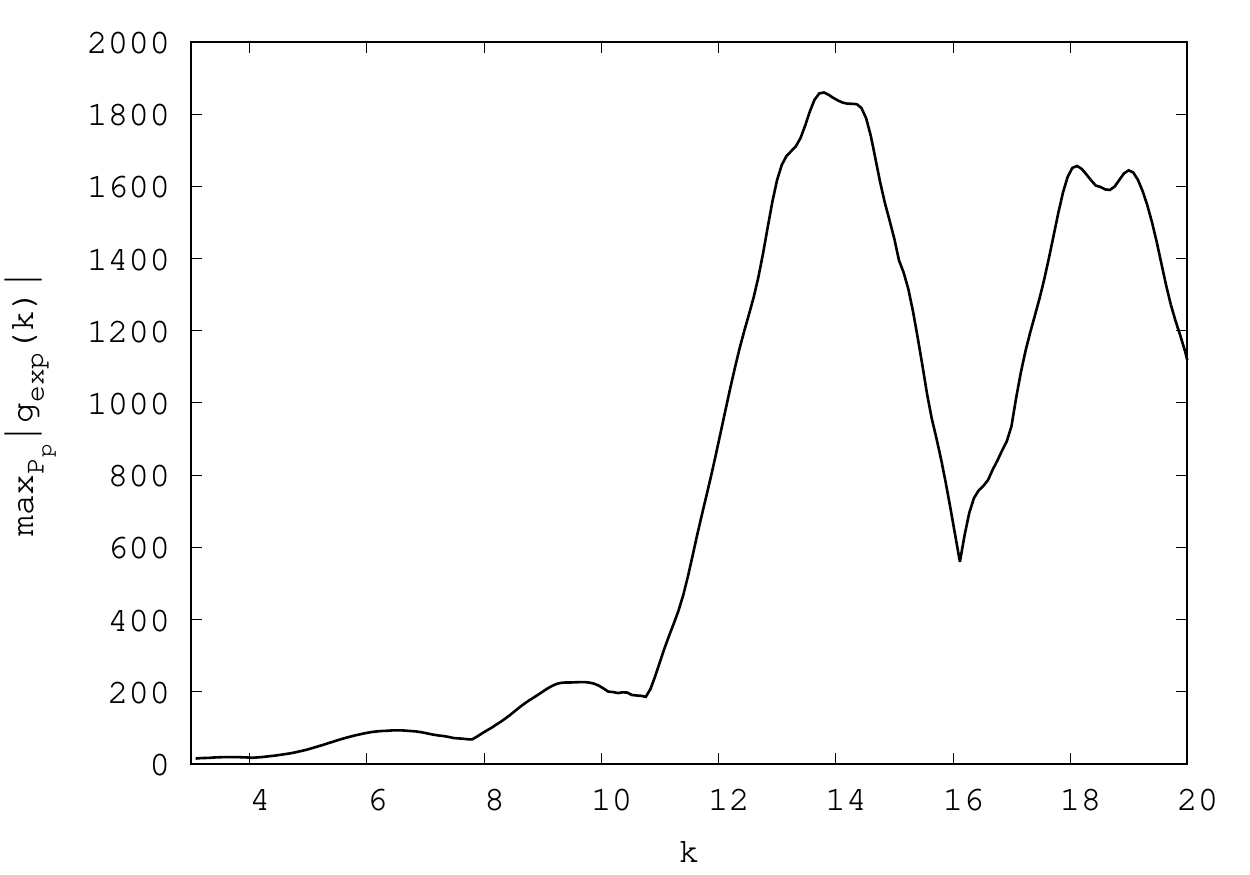}} %
\subfloat[\label{fig 7a}]{\includegraphics[width=0.4\textwidth]{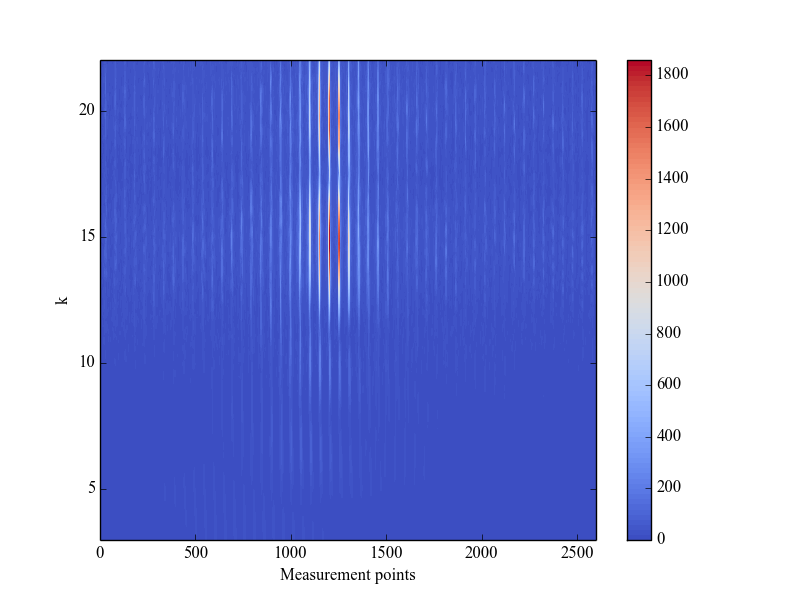}}
\end{center}
\caption{\label{fig:g_max_obj1} The graph of the function $s\left( k\right) $ defined in (%
\ref{1010}) for the reference target number 1. (b) The absolute value of the
propagated experimental data. The horizontal line means grid points of the
rectangle $S_{p}$ on the propagated plane $P_{p},$ see (\ref%
{102}). The vertical axis denotes the wave numbers. We observe on (b)
that the absolute values of the propagated data for $k\in \lbrack 13.5,14]$
are near the maximal value of the function $|\tilde{g}(x,k)|$
\ and these data strongly focus in a subdomain of the rectangle $S_{p}$.}
\end{figure}

\subsection{Complementing the backscatter data}

\label{sec:5.2}

Our backscatter experimental data are measured only on the rectangle $S_{m}$
of the measurement plane $P_{m}=\left\{ x_{3}=-8\right\} $ and then are
propagated to the plane $P_{p}=\left\{ x_{3}=-0.75\right\} $. Hence, we
denote 
\[
\Gamma =\partial \Omega \cap \left\{ x_{3}=-0.75\right\} . 
\]%
In fact, $\Gamma =S_{p}\subset P_{p}$. On the other hand, our algorithm of
section 3.3 works only with the case when the data are given at the entire
boundary $\partial \Omega .$ Thus, using (\ref{1000}), we complement the
missing data on the rest of the boundary $\partial \Omega $ for $k\in \left[
6,6.5\right] $ by the data in air $c(x)=1$ 
\begin{equation}
\widehat{g}(x,k)=\left\{ 
\begin{array}{c}
e^{ikx_{3}}+d\left( k\right) \tilde{g}_{exp}(x,k+\underline{k}^{\prime
}-6),x\in \Gamma , \\ 
e^{ikx_{3}},x\in \partial \Omega \diagdown \Gamma .%
\end{array}%
\right.  \label{25}
\end{equation}%
This formula can be justified assuming that the targets are located far from 
$\partial \Omega \setminus \Gamma $ and the influence of waves scattered
from $\partial \Omega \setminus \Gamma $ is small. We note that in the above
cited references \cite{IPexp1,KLN,Liem,TBKF1,TBKF2} the backscatter data
were complemented the same way as in (\ref{25}). The $k-$derivative $%
\partial _{k}\tilde{g}_{exp}$ was found using the finite difference method
with the step size $h=0.05$, which is the same as in Remark 4.1. Even though
the differentiation of noisy data is an ill-posed problem, we have not
observed any instabilities. This is because of the efficiency of the data
propagation and the choice of the optimal interval of wave numbers. Note
that a similar differentiation procedure, although with respect to the
parameter of the Laplace transform, was applied in \cite{BK,BK1} to
computationally simulated data with noise and to experimental data in \cite%
{IPexp1,TBKF1,TBKF2}. In the case of the differentiation with respect to $k$%
, it was also applied in \cite{KLN,Liem}. It was stable in all cases. A
further study of the differentiation topic is outside of the scope of this
paper.

\subsection{The first tail function $V_{0}\left( x\right) $}

\label{sec:5.3}

Since we set $\Delta V_{0}=0$ (section 3.4), then it is clear from (\ref{12}%
) (at $n=1$) that we use only the gradient $\nabla V_{0}$ of the first tail
function $V_{0}$. Thus, to avoid the error associated with the numerical
differentiation of the first tail function $V_{0}(x)$, we numerically solve
the following problem instead of the problem (\ref{180}), (\ref{190}) 
\[
\Delta (\partial _{x_{j}}V_{0})=0,\quad x\in \Omega . 
\]%
\begin{equation}
\partial _{x_{j}}V_{0}=\frac{\partial _{x_{j}}u(x,\overline{k})}{u(x,%
\overline{k})},\quad x\in \Gamma ,\quad j=1,2,3,  \label{103}
\end{equation}%
\begin{equation}
\partial _{x_{j}}V_{0}=0,\quad j=1,2,\quad x\in \partial \Omega \setminus
\Gamma ,  \label{104}
\end{equation}%
\begin{equation}
\partial _{x_{3}}V_{0}=i\overline{k},\quad x\in \left\{ x_{3}=4.25\right\}
\cap \partial \Omega .  \label{105}
\end{equation}%
We now comment about boundary conditions (\ref{103})-(\ref{105}). Conditions
(\ref{103}) follow from (\ref{10}), (\ref{200}), (\ref{14}) and (\ref{190}).
Conditions (\ref{104}) and (\ref{105}) for $x\in \left\{ x_{3}=4.25\right\}
\cap \partial \Omega $ follow from our assumption in the second line of (\ref%
{14}). In addition, conditions (\ref{104}) are approximate ones for those
parts of the rest of $\partial \Omega \setminus \Gamma $ which are
orthogonal to the $x_{j}-$axis. In the latter case they basically mean that
the function $V_{0}$ does not change in the direction $x_{j}$ near of those
parts of the boundary. As to condition (\ref{103}), it implies that we need
the boundary data for functions $\partial _{x_{j}}u(x,\overline{k})$, $x\in
\Gamma ,$ $j=1,2,3$. The derivatives with respect to $x_{1}$ and $x_{2}$ are
found by differentiating $g(x,\overline{k})$ which is defined in (\ref{1000}%
). The derivative with respect to $x_{3}$ for $x\in \Gamma $ is calculated
as 
\[
\partial _{x_{3}}u(x,\overline{k})=d(\overline{k})\partial _{x_{3}}\tilde{g}%
(x,\overline{k})+i\overline{k}\exp (i\overline{k}x_{3}), 
\]%
where $\partial _{x_{3}}\tilde{g}(x,\overline{k})$ is the derivative of the
propagated experimental data $\tilde{g}(x,\overline{k})$ with respect to $%
x_{3}$. The function $\partial _{x_{3}}\tilde{g}(x,\overline{k})$ was found
by propagating the Fourier transform of experimental data to the two nearby
planes $P_{p}=\{x_{3}=-0.75\}$ and $P_{p_{\varepsilon
}}=\{x_{3}=-0.75+\varepsilon \}$, subtracting the results from each other
and dividing by $\varepsilon $, where we took $\varepsilon =0.1$.
Derivatives $\partial _{x_{j}}g(x,\overline{k}),j=1,2$ for $x\in \Gamma $
were found using FreeFEM++ and Matlab. The preprocessed data looked very
smooth and we have not observed any instabilities.

\subsection{Truncation}

\label{sec:5.4}

As it is clear from results of data propagation (Figures 5b,d and 6b), $%
x_{1},x_{2}$ coordinates of targets can be estimated from images on the
propagated plane $P_{p}$. We have observed for all targets of our study that
the absolute value of the Fourier transform of the propagated experimental
data $\tilde{g}(x,k)$ has a positive peak (see, e.g. Fig.~\ref{fig:g}),
which corresponds to the $x_{1},x_{2}$ location of a target. Thus, we define
the set $\Gamma _{t}\subset S_{p}$ as 
\[
\Gamma _{t}=\left\{ (x_{1},x_{2}):|\tilde{g}(x,\overline{k})|>0.7\max |%
\tilde{g}(x,\overline{k})|,(x_{1},x_{2})\in S_{p}\right\} . 
\]%
Here, the truncation value 0.7 is chosen by trial and error as an optimal
number.

Taking into account our desired application to detection and identification
of explosives, we look for our targets only in the interval $(-0.75,1)$ in
the $x_{3}-$direction. This means that we allow the maximal linear size of
the target in the $x_{3}$-direction to be 17.5 cm, which is close to the
upper bound of about 15 cm for antipersonnel mines and IEDs. Recall that the
function $c_{n,i}(x)$ is calculated by formula (\ref{13}) in our algorithm.
However, the right hand side of this formula might be complex valued. Hence,
we need to postprocess (\ref{13}). So, we define the postprocessed function $%
c_{n,i}^{\left( p\right) }(x)$ as: 
\begin{equation}
c_{n,i}^{\left( p\right) }(x)=\left\{ 
\begin{array}{cc}
\left\vert c_{n,i}(x)\right\vert , & \quad x\in \Gamma _{t}\times
(-0.75,1.0), \\ 
1, & \quad \mbox{elsewhere},%
\end{array}%
\right.  \label{106}
\end{equation}%
where $c_{n,i}(x)$ is the function calculated by (\ref{13}). Next, the
coefficient $c_{n,i}^{\left( p\right) }$ is smoothed by the \texttt{smooth3}
function in Matlab. That smoothed coefficient is used then to solve the
Lippmann-Schwinger equation in our algorithm.

In the case of targets with two inclusions, their $x_{1},x_{2}$ coordinates
are seen on the propagated plane $P_{p}$. However, the signal from one
inclusion might be stronger than from the other, see Fig. \ref{fig:g2}. This
might mean that either one inclusion has a larger dielectric constant than
the second one, or perhaps they have different volumes, or both. Thus, it is
hard to reconstruct both the \textquotedblleft stronger" and the
\textquotedblleft weaker" inclusions simultaneously. Therefore, we treat
these inclusions as two separate targets. For example, in Fig. \ref{fig:g2}
we consider the median straight line on the plane $P_{p}$ between images of
these two inclusions on $P_{p}$. Next, to image the top inclusion, we set
the propagated data below this line to be zero $\tilde{g}(x,k):=0$. Next, we
apply the above procedure, which gives us the image of the top inclusion.
Having the image of the top inclusion, we then set $\tilde{g}(x,k):=0$ for
points $x\in P_{p}$ located above that median line and repeat. Finally we
combine two images in one. The reconstruction results for both inclusions
look good and their $x_{3}-$locations are the same as it should be, see
Figures 11 and 12.

\section{Reconstruction results}

\label{sec:6}

The preprocessed experimental data (\ref{25}), as well as the gradient of
the first tail function, which was calculated as in section 5.3, were used
as the input for our globally convergent algorithm of subsection 3.3. We
present in this section results of the reconstructions by this algorithm.

We have considered ten (10) data sets. In seven cases we had one target, and
we had two targets in the rest of three data sets. Six targets were
dielectrics, three were metallic ones and one target was a mixture of a wood
and metal. Although the dielectric constant of the metal is not
well-defined, it was established in an earlier work \cite{KSNF1} that, in
the case of metallic targets, one can assume that they have large the
so-called \emph{appearing} dielectric constants, 
\begin{equation}
c\left( \mbox{metal}\right) \in \left[ 10,30\right] .  \label{1}
\end{equation}%
A description of targets and their linear sizes is given in Table \ref%
{tab:targers}. Recall that the target number 1 was chosen as the reference
object for data calibration.

\begin{table}[tbp]
\begin{center}
\begin{tabular}{|c|l|c|}
\hline
Target & Description & Size in cm, \\ 
number &  & ($x_1 \times x_2 \times x_3$) \\ \hline
1 & A piece of oak & $4.1 \times 8.2 \times 4.1 $ \\ \hline
2 & Metallic ball & $10.3 \times 10.3 \times 10.3$ \\ \hline
3 & Metallic cylinder & $5.3 \times 11.1 \times 5.3$ \\ \hline
4 & Wooden object & $6.0 \times 11.3 \times 4.0 $ \\ \hline
5 & Wooden doll with air inside & $8.7 \times 11.5 \times 8.7$ \\ \hline
6 & Wooden doll with metal inside & $8.7 \times 11.5 \times 8.7$ \\ \hline
7 & Wooden doll filled with sand & $8.7 \times 11.5 \times 8.7$ \\ \hline
8 & Two metallic objects & $5.5 \times 10.1 \times 2.5$ \\ 
&  & $5.9 \times 4.1 \times 4.1$ \\ \hline
9 & Two wooden objects aligned vertically & $9.6 \times 5.8 \times 5.3$ \\ 
&  & $11.3 \times 6.0 \times 3.3$ \\ \hline
10 & Two wooden objects aligned horizontally & $4.1 \times 8.15 \times 4.1$
\\ 
&  & $5.7 \times 9.7 \times 5.7$ \\ \hline
\end{tabular}%
\end{center}
\caption{Targets and their linear sizes in centimeters}
\label{tab:targers}
\end{table}

We choose the final reconstructed coefficient as follows. First, we
calculate the following relative errors%
\begin{equation}
\varepsilon _{n,i}=\left\{ 
\begin{array}{c}
\frac{\Vert c_{n,i}-c_{n,i-1}\Vert _{L_{2}(\Omega )}}{\Vert c_{n,i-1}\Vert
_{L_{2}(\Omega )}},\mbox{ for }i=2,...m,\,1\leq n\leq N=10, \\ 
\frac{\Vert c_{n+1,1}-c_{n}\Vert _{L_{2}(\Omega )}}{\Vert c_{n}\Vert
_{L_{2}(\Omega )}},\mbox{ for }i=1,\,1\leq n\leq N=9.%
\end{array}%
\right.  \label{26}
\end{equation}%
Suppose that the minimal error out of numbers (\ref{26}) is achieved at $%
n=n_{0},i=i_{0}.$ Then we set the computed coefficient $c_{comp}\left(
x\right) =c_{n_{0}i_{0}}\left( x\right) .$ If the minimal error is achieved
for several pairs $\left( n,i\right) ,$ then we choose the pair $\left(
n_{0},i_{0}\right) $ which was the first one in (\ref{26}). Note that the
number of inner iterations in our computations is $m=3$. This number was
chosen by trial and error.

To verify the accuracy of our computations, we have a posteriori measured
directly refractive indices $n=\sqrt{c}$ of our six dielectric targets.
Since the latter measurements were in time domain, the dependencies of
measured indices on the frequency were not known, but we believe that
approximating it by a constant in the small intervals of wave numbers we
work with is reasonable. In table \ref{tab:results_nonmetal} we present the
measured and computed refractive indices of the dielectric targets together
with the errors in their measurements and computational errors. The computed
value of refractive index is defined as $n_{comp}=\sqrt{\max_{\overline{%
\Omega }}c_{comp}(x)}$ and its error is computed as 
\[
\varepsilon _{comp}=\frac{\left\vert n_{true}-n_{comp}\right\vert }{n_{true}}%
100\%. 
\]

Table \ref{tab:results_nonmetal} demonstrates that our computational results
for refractive indices of dielectric targets are highly accurate. Indeed,
the computational errors are smaller than the corresponding measurement
errors. In fact, the error for the reference target 1 is the smallest one,
which indicates that the experimental data for this target were correctly
scaled using the calibration factor (\ref{eq:calib_factor}). We also observe
that the computational errors for the second object of data sets 9 and 10
are higher than the errors for the first object in the same data set. This
can be explained by the fact that the second objects in these data sets have
smaller refractive indices. Hence, their the reflected signals are weaker
than ones of first objects, which, in turn leads to a higher computational
error.

Recall that we associate the so-called appearing dielectric constant with
metals, see (\ref{1}). Thus, we model metallic targets as the dielectric
targets with high values of the dielectric constant $c\in \lbrack 10,30]$.
In table \ref{tab:results_metal} we present the computed appearing
dielectric constants $c_{c}=\max_{\overline{\Omega }}c_{comp}(x)$ for our
metallic targets. We see that the computed appearing dielectric constant $%
c_{c}=6.68$ of the target number 2 is less than 10. This is because the
signal reflected by this target is not strong enough. Indeed, for this
target, the maximal absolute value of the preprocessed time dependent signal
after Step 3 was about 800. On the other hand, the analogous value for
targets number 3 and 8 were about 1800 and 3200, respectively. 

The maximal absolute value of the preprocessed time dependent signal after
Step 3 for the target number 6 was about 1100, which is also significantly
less than this value for targets number 3 and 8. The computed appearing
dielectric constant 9.73 of the target number 6 is close to the lower limit
of 10. It is natural, however, that it is lower than 10, since a piece of
metal in this target was covered by the surface of that wooden doll. The
latter, in turn has caused the decrease of that maximal absolute value.

\begin{figure}[tbp]
\centering
\subfloat[\label{fig:comp_obj1}]{\includegraphics[width=0.46%
\textwidth]{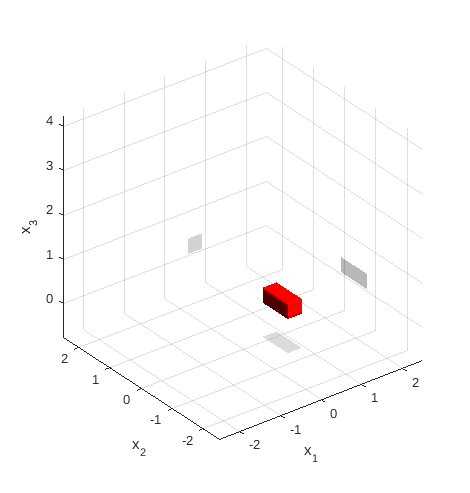}} \subfloat[\label{fig:true_obj1}]{%
\includegraphics[width=0.52\textwidth]{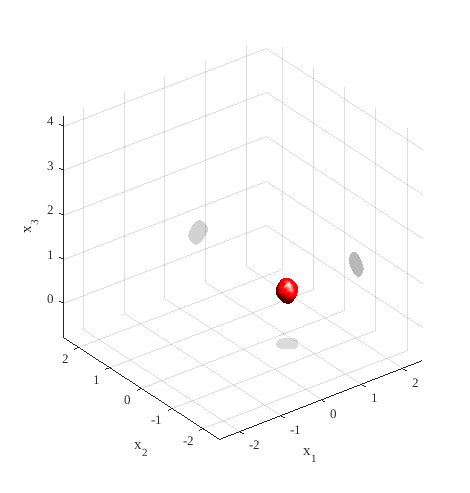}}
\caption{Reconstruction result for the target number 1: true image (a) and
computed image (b).}
\label{fig:res_obj1}
\end{figure}

\begin{figure}[tbp]
\centering
\subfloat[\label{fig:comp_obj3}]{\includegraphics[width=0.49%
\textwidth]{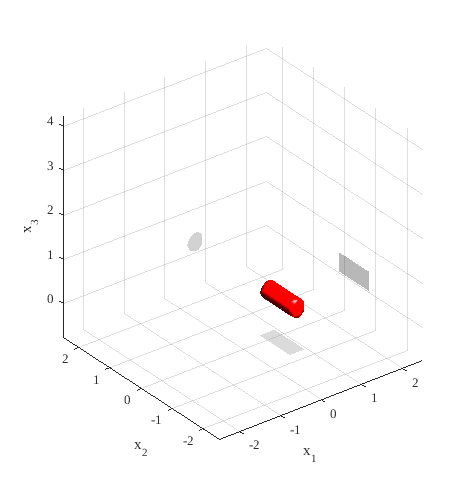}} \subfloat[\label{fig:true_obj3}]{%
\includegraphics[width=0.49\textwidth]{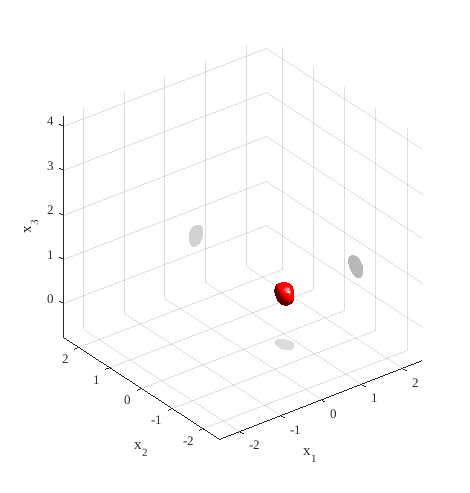}}
\caption{Reconstruction result for the target number 3: true image (a) and
computed image (b).}
\label{fig:res_obj3}
\end{figure}

\begin{figure}[tbp]
\centering
\subfloat[\label{fig:comp_obj5}]{\includegraphics[width=0.49%
\textwidth]{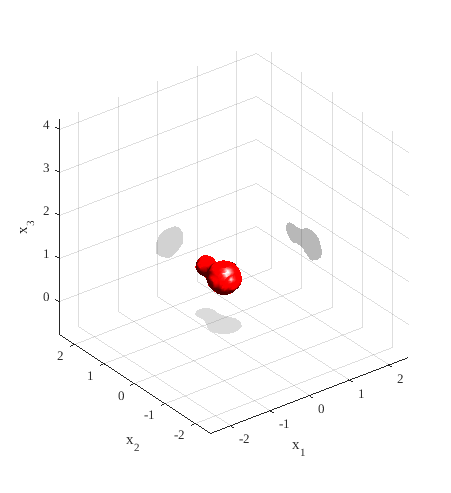}} \subfloat[\label{fig:true_obj5}]{%
\includegraphics[width=0.49\textwidth]{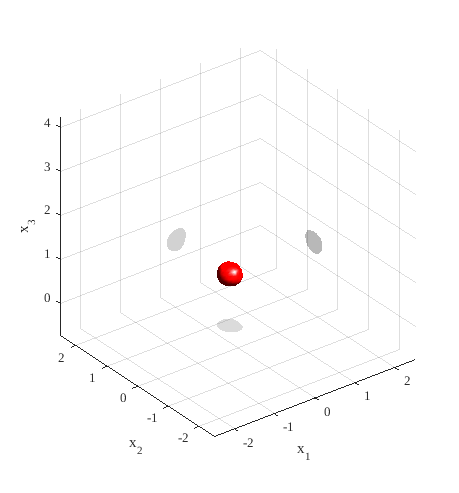}}
\caption{Reconstruction result for the target number 5: true image (a) and
computed image (b).}
\label{fig:res_obj5}
\end{figure}

\begin{figure}[tbp]
\centering
\subfloat[\label{fig:comp_obj8}]{\includegraphics[width=0.49%
\textwidth]{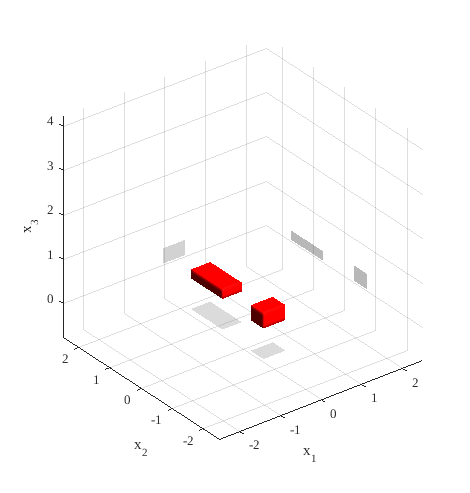}} \subfloat[\label{fig:true_obj8}]{%
\includegraphics[width=0.49\textwidth]{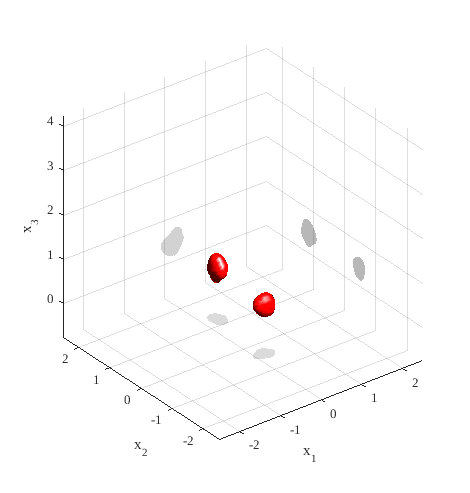}}
\caption{Reconstruction result for the target number 8: true image (a) and
computed image (b).}
\label{fig:res_obj8}
\end{figure}

\begin{figure}[tbp]
\centering
\subfloat[\label{fig:comp_obj10}]{\includegraphics[width=0.49%
\textwidth]{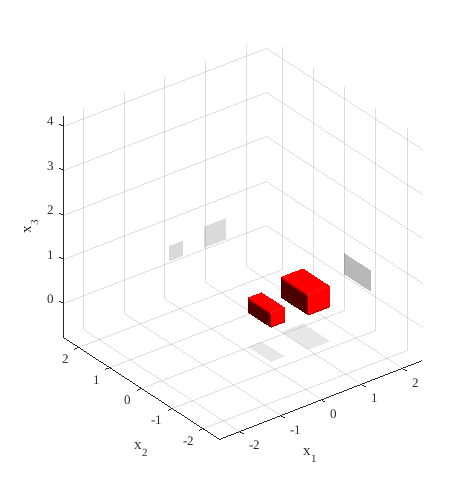}} \subfloat[\label{fig:true_obj10}]{%
\includegraphics[width=0.49\textwidth]{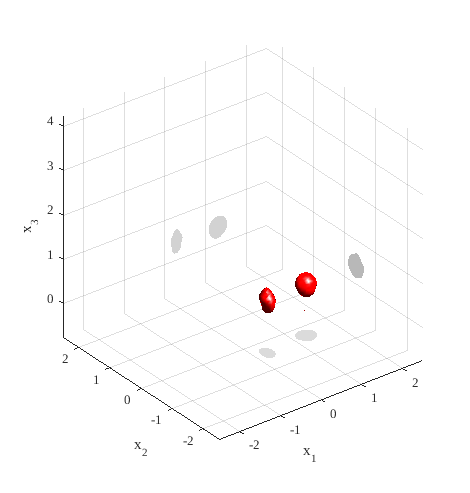}}
\caption{Reconstruction result for the target number 10: true image (a) and
computed image (b). Note that the computed image corresponding to the
inclusion of a larger size on (a) is indeed larger on (b).}
\label{fig:res_obj10}
\end{figure}

Figures \ref{fig:res_obj1} -- \ref{fig:res_obj10} illustrate the true and
computed images for the targets number 1, 3, 5, 8, 10. We use \textit{%
isosurfaces} of Matlab to draw these figures. Reconstructed images for the
rest of five targets with numbers 2, 4, 6, 7, 9 are similar. We see that the
positions of objects are imaged correctly. Shapes are not well
reconstructed. Still, in the case of two objects of target number 10 (Fig. %
\ref{fig:comp_obj10}), one object is larger than another one, and we observe
the same on computed objects (Fig. \ref{fig:true_obj10}), which is accurate.
In the case of the wooden doll with air inside of the object number 5 (Fig. %
\ref{fig:comp_obj5}), only its larger part is reconstructed, which is
natural.

\begin{table}[tbp]
\begin{center}
\begin{tabular}{|c|c|c|}
\hline
Target & \multicolumn{2}{c|}{Refractive index} \\ \cline{2-3}
number & Measured value $n_{true}$ (error, \%) & Computed value $n_{comp}$
(error, \%) \\ \hline
1 & 2.11 $\quad$ (19) & 2.08 $\quad$ (1.39) \\ \hline
4 & 2.14 $\quad$ (28) & 2.54 $\quad$ (18.8) \\ \hline
5 & 1.89 $\quad$ (30) & 2.06 $\quad$ (8.99) \\ \hline
7 & 2.10 $\quad$ (26) & 2.28 $\quad$ (8.38) \\ \hline
9 & 2.14 $\quad$ (18) & 2.25 $\quad$ (5.21) \\ 
& 1.84 $\quad$ (18) & 2.03 $\quad$ (10.3) \\ \hline
10 & 2.14 $\quad$ (18) & 2.28 $\quad$ (6.48) \\ 
& 1.84 $\quad$ (18) & 2.08 $\quad$ (13.3) \\ \hline
\end{tabular}%
\end{center}
\caption{Measured and computed refractive indices of dielectric targets.}
\label{tab:results_nonmetal}
\end{table}

\begin{table}[tbp]
\begin{center}
\begin{tabular}{|c|c|}
\hline
Target & $\quad$ Computed appearing dielectric constant $c_{comp}\quad $ \\ 
number &  \\ \hline
2 & 6.68 \\ \hline
3 & 11.49 \\ \hline
6 & 9.73 \\ \hline
8 & 18.47 \\ 
& 18.73 \\ \hline
\end{tabular}%
\end{center}
\caption{Computed appearing dielectric constants of metallic targets}
\label{tab:results_metal}
\end{table}

\section{Discussion}

\label{sec:7}

In this paper, we have tested the performance of the GCM of \cite{KLN} on
time dependent experimental data. Since the technique of \cite{KLN} works
for multi-frequency data, we have applied the Fourier transform to
preprocessed time dependent data. The main difficulty working with the data
of both the current paper and \cite{Liem} was a \emph{huge discrepancy}
between these data and computationally simulated ones. Standard denoising
techniques do not work for this case due to a richer content of the real
data, compared with the computationally simulated ones. Thus, we have
developed here a heuristic data preprocessing procedure. The resulting data
look somehow similar to the computationally simulated ones. These
preprocessed data are used then as the input for the GCM.

Our data preprocessing procedure consists of two stages. On the first stage,
the time dependent data are preprocessed.\ On the second stage, the Fourier
transform is applied to the data preprocessed on the first stage, and then
the resulting data are preprocessed again in the frequency domain. Even
though the first stage is similar with \cite{TBKF1} and the second stage is
similar with \cite{Liem}, there are two important differences with these two
references: one difference in each stage. The difference in the first stage
is caused by the difference between the Fourier transform used in the
current publication and the Laplace transform used in \cite{TBKF1}. The
difference in the second stage is in the data shift (\ref{22}), which was
not used in \cite{Liem}. This data shift is caused by the difference between
our Fourier transformed data and the multi-frequency experimental data of 
\cite{Liem}.

Some elements of our heuristic data preprocessing procedure might seem to be
of a concern, such as, e.g. formula (\ref{22}). However, we believe that the
ultimate judgement about the usefulness of this procedure as a whole should
be drawn on the basis of obtained reconstruction results. As to the latter,
it is clear from table \ref{tab:results_nonmetal} that numerical results of
reconstructions of refractive indices of dielectric targets are very
accurate. Furthermore, they are even significantly more accurate than
results of direct measurements of $n$. Recall that refractive indices $n$,
rather than dielectric constants $c=n^{2}$ were a posteriori directly
measured in our experiments. We also recall that it is hard for optimization
methods for CIPs to accurately calculate values of unknown coefficients
inside the targets.

Table \ref{tab:results_metal} demonstrates that appearing dielectric
constants $c$ of metallic targets number 3 and 8 are also in the required
range of $c\in \left[ 10,30\right] ,$ see (\ref{1}). Although the computed
value $c_{c}=6.68$ for the metallic target number 2 is less than 10, it is
still quite high, as compared with dielectric constants of non metallic
targets in table \ref{tab:results_nonmetal}, all of which do not exceed 2.6.
As it was explained in section 6, the reason why the computed appearing
dielectric constant of target number 2 is significantly less than 10 is that
the preprocessed time dependent signal after Step 3 for this case was
significantly weaker than the one for the targets number 3 and 8.

The fact that the computed appearing dielectric constant $c_{c}=9.73$ of the
target number 6 is less than 10 has a natural explanation. Indeed, this
target is a piece of metal hidden inside of an otherwise empty wooden doll,
see Table \ref{tab:targers}. Therefore, the surface of this doll works as a
\textquotedblleft mask" for that metal. This observation seems to be of an
interest for our desired application to imaging of mines and IEDs.

Locations of all targets are imaged accurately. Although shapes are computed
inaccurately, in the case of two objects of Figure. \ref{fig:comp_obj10},
the smaller object is computed as the smaller one, see Figure \ref%
{fig:true_obj10}.

The main advantage of our GCM is that, given a reasonable mathematical
assumption, there is a rigorous guarantee that it provides some points in a
small neighborhood of the exact coefficient \cite{KLN}. This is an important
difference compared with conventional optimization methods. A reasonable
mathematical assumption seems to be unavoidable due to the well known fact
that the development of globally convergent numerical methods for CIPs, all
of which are highly nonlinear, is a substantially challenging task,
especially in the case of single measurement data. Our reasonable
mathematical assumption is (\ref{150}) and (\ref{16}), which amounts to
ignoring the small term $O(1/k)$ in (\ref{14}). That assumption is used only
once: for defining the first tail function $V_{0}.$ However, it is not used
on follow up iterations of our numerical method.

One of conditions of the \ global convergence theorem for our GCM \cite{KLN}
is that the interval $\left[ \underline{k},\overline{k}\right] $ of wave
numbers should be sufficiently small, so as the discretization step size $h$
with respect to $k$. Interestingly, this is exactly the case of our data, so
as of the data of \cite{Liem}. Indeed, consider the relative length of this
interval $a=\left( \overline{k}-\underline{k}\right) /\underline{k}.$ Then
in our case $a=0.5/6=0.083$ after the data shift (\ref{22}) (see section
4.2.7). Hence, in our case $h=0.0083=0.083/10.$ Both these numbers are
sufficiently small. In addition, the analysis of section 4.2.7 indicates
that it is likely meaningless to work on a larger interval of wave numbers.
This provides a sort of a justification of that smallness condition from the
Physics standpoint.

\begin{center}
\textbf{Acknowledgements}
\end{center}

This work was supported by the US Army Research Laboratory and US Army
Research Office grant W911NF-15-1-0233 as well as by the Office of Naval
Research grant N00014-15-1-2330.


\end{document}